\documentclass[11pt]{article}
\usepackage{amsfonts}
\usepackage{color}
\newtheorem{theo}{Theorem}[section]
\newtheorem{lem}[theo]{Lemma}

\newcommand{\mysection}[1]{\section{#1} \setcounter{equation}{0}}
\newcommand{\proof}{{\sc Proof.} \quad}
\newcommand{\proofc}{{\sc Proof} \ }
\newcommand{\be}{\begin{equation} \label}
\newcommand{\ee}{\end{equation}}
\newcommand{\bea}{\begin{eqnarray}\label}
\newcommand{\eea}{\end{eqnarray}}
\newcommand{\bas}{\begin{eqnarray*}}
\newcommand{\eas}{\end{eqnarray*}}
\newcommand{\bit}{\begin{itemize}}
\newcommand{\eit}{\end{itemize}}
\newcommand{\qed}{\hfill$\Box$ \vskip.2cm}
\newcommand{\nn}{\nonumber}
\newcommand{\R}{\mathbb{R}}
\newcommand{\N}{\mathbb{N}}
\newcommand{\pO}{\partial\Omega}

\newcommand{\eps}{\varepsilon}

\newcommand{\wto}{\rightharpoonup}
\newcommand{\wsto}{\stackrel{\star}{\rightharpoonup}}
\newcommand{\hra}{\hookrightarrow}
\newcommand{\io}{\int_\Omega}

\newcommand{\mult}{\otimes}
\newcommand{\bom}{\overline{\Omega}}

\newcommand{\abs}{\\[5pt]}

\newcommand{\neps}{n_\eps}
\newcommand{\ceps}{c_\eps}
\newcommand{\ueps}{u_\eps}
\newcommand{\Peps}{P_\eps}
\newcommand{\yeps}{y_\eps}
\newcommand{\zeps}{z_\eps}
\newcommand{\heps}{h_\eps}
\newcommand{\nee}{n_{\eps\eta}}
\newcommand{\cee}{c_{\eps\eta}}
\newcommand{\uee}{u_{\eps\eta}}
\newcommand{\Pee}{P_{\eps\eta}}
\newcommand{\proj}{{\mathcal{P}}}
\newcommand{\hatq}{\widehat{q}}
\newcommand{\ns}{{\mathcal{N}}}
\newcommand{\zho}{z_h^{(1)}}
\newcommand{\zht}{z_h^{(2)}}
\begin{document}
\enlargethispage{10mm}
\title{The fast signal diffusion limit in Keller-Segel(-fluid) systems}
\author{
Yulan Wang\footnote{wangyulan-math@163.com}\\
{\small  School of Science,
Xihua University, Chengdu 610039, China}\\
\and
Michael Winkler\footnote{michael.winkler@math.uni-paderborn.de}\\
{\small Institut f\"ur Mathematik, Universit\"at Paderborn,  }\\
{\small 33098 Paderborn, Germany}
\and
Zhaoyin Xiang\footnote{zxiang@uestc.edu.cn}\\
{\small School of Mathematical Sciences, University of Electronic Science and Technology of China,}\\
{\small Chengdu 611731, China}
}
\date{}
\maketitle
\begin{abstract}
\noindent
This paper deals with convergence of solutions to a class of parabolic Keller-Segel systems, possibly coupled to
the (Navier-)Stokes equations in the framework of the full model
\begin{eqnarray*}
	\left\{ \begin{array}{lcl}
	\, \, \partial_t \neps + \ueps \cdot \nabla\neps
	&=& \Delta \neps - \nabla \cdot \Big( \neps S(x,\neps,\ceps)\cdot\nabla \ceps\Big) + f(x,\neps,\ceps),  \\[1mm]
	\eps \partial_t \ceps + \ueps\cdot\nabla \ceps
	&=& \Delta \ceps - \ceps + \neps, \\[1mm]
	\,\,\partial_t \ueps + \kappa (\ueps\cdot\nabla)\ueps
	&=&  \Delta \ueps + \nabla \Peps + \neps \nabla\phi,
	\qquad \nabla\cdot \ueps=0
	\end{array} \right.
\end{eqnarray*}
to solutions of the parabolic-elliptic counterpart formally obtained on taking $\eps\searrow 0$.
In smoothly bounded physical domains $\Omega\subset {\mathbb R}^{N}$ with $N\ge 1$, and under
appropriate assumptions on the model ingredients, we shall first derive a general result which asserts certain
strong and pointwise convergence properties whenever asserting that
supposedly present bounds on
$\nabla\ceps$ and $\ueps$ are bounded in $L^\lambda((0,T);L^q(\Omega))$ and in
$L^\infty((0,T);L^r(\Omega))$, respectively, for some $\lambda\in (2,\infty]$, $q>N$ and $r>\max\{2,N\}$ such that
$\frac{1}{\lambda}+\frac{N}{2q}<\frac{1}{2}$.
To our best knowledge, this seems to be the first rigorous mathematical result on a fast signal diffusion
limit in a chemotaxis-fluid system.\abs
This general result will thereafter be concretized in the context of two examples:
Firstly, for an unforced Keller-Segel-Navier-Stokes system we shall
establish a statement on global classical solutions under suitable smallness conditions on the initial data,
and show that these solutions approach a global classical solution to the respective parabolic-elliptic simplification.\abs
We shall secondly derive a corresponding convergence property for arbitrary solutions to
fluid-free Keller-Segel systems with logistic source terms, which in spatially one-dimensional settings turn out
to allow for a priori estimates compatible with our general theory.
Building on the latter in conjunction with a known result on emergence of large densities in the associated
parabolic-elliptic limit system,
we will finally discover some quasi-blowup phenomenon for the fully parabolic Keller-Segel system with
logistic source and suitably small parameter $\eps>0$.\abs
\end{abstract}
\newpage
\mysection{Introduction}\label{intro}
{\bf The Keller-Segel system and its parabolic-elliptic simplification.} \quad
To describe chemotactic aggregation of cellular slime molds which move towards relatively high concentrations of a chemical secreted by the amoebae themselves,
Keller and Segel \cite{keller_segel} proposed cross-diffusive parabolic systems of the form
\begin{eqnarray*}
	\left\{ \begin{array}{lcl}
	n_t &=& d_1 \Delta n - a_1 \nabla \cdot (n\nabla c), \\[1mm]
	c_t &=& d_2\Delta c-a_2c+a_3n,
		\end{array} \right.
\end{eqnarray*}
where the unknown functions $n=n(x,t)$ and $c=c(x,t)$ denote the cell density and the concentration of the
chemical substance at place $x$ and time $t$, respectively, and where $d_1,\, d_2,\,a_1, \, a_2,\, a_3$ are
positive numbers. By substituting
\bas
	\frac{a_1}{d_1}=S,\quad \frac{d_1}{d_2}=\varepsilon, \quad \frac{a_2}{d_2}=\gamma \quad
	\mbox{and} \quad
	\frac{a_3}{d_2}=\alpha,
\eas
and replacing $d_1 t$ with $t$, from this we obtain the system
\be{KS}
	\left\{ \begin{array}{lcl}
	\,\, n_t &=& \Delta n - \nabla \cdot (n S\nabla c),   \\[1mm]
	\varepsilon c_t &=& \Delta c-\gamma c+\alpha n,
	\end{array} \right.
\ee
which in the limit $\eps\searrow 0$ formally approaches the corresponding parabolic-elliptic system, with the second
identity therein replaced with the inhomogeneous Helmholtz equation $-\Delta c+ \gamma c=\alpha n$.\abs
As is well-known from quite a large literature, with regard to technical purposes the latter simplification
goes along with substantial advantages for mathematical analysis, in summary leading to much a deeper knowledge
on parabolic-elliptic Keller-Segel systems than currently available for their fully parabolic relatives.
Examples already include the mere detection of exploding solutions, typifying the probably most characteristic
effect of the considered cross-diffusive interaction, which in fact could be accomplished for parabolic-elliptic
systems already rather early (\cite{JL}, \cite{nagai1995}, \cite{nagai2001}, \cite{biler}), while for the
full system (\ref{KS}) with positive $\eps$, corresponding results on generic blow-up, thus going
beyond particular examples (\cite{herrero_velazquez}), seem to require significantly stronger efforts and hence
have been achieved only a few years ago (\cite{win_JMPA}, \cite{mizoguchi_win}).
Likewise, while considerable qualitative knowledge on the respective blow-up mechanisms has been collected
for parabolic-elliptic systems (see e.g.~\cite{NSS2000}, \cite{naito_suzuki2008}, \cite{senba2007},
\cite{suzuki_JMPA2013}, \cite{suzuki_book2005}, 	
\cite{souplet_win}), only little information
seems available for general blow-up solutions to (\ref{KS}) when $\eps>0$
(\cite{NSS2000}, \cite{mizoguchi_souplet}, \cite{win_ct_profile}).\abs
More generally, by providing accessibility to numerous tools, especially from the analysis of scalar parabolic problems,
resorting to parabolic-elliptic simplifications has		
made it possible to reveal further qualitative properties of Keller-Segel-type systems,
inter alia also in the framework of global solutions (\cite{kavallaris_souplet}, \cite{kang_stevens},
\cite{nadin_perthame_ryzhik}, \cite{win_ct_critmass}), and partially even including couplings to additional
quantities such as fluid flows or haptotactic attractants (\cite{kiselev_ryzhik_CPDE}, \cite{kiselev_ARMA2016},
\cite{taowin_hapto}).\abs
{\bf Problem setting and main objectives.} \quad
In line with the above, it seems natural to seek for some appropriate control of the error made when approximating
a fully parabolic system of Keller-Segel type by its parabolic-elliptic simplification, especially in cases
when the considered biological situation is such that the respective signal diffuses much faster than individuals
in the cell population, in the context of (\ref{KS}) thus meaning that $\eps>0$ is small.
Indeed, even in the context of the classical system (\ref{KS}) already the question concerning mere convergence
of solutions as $\eps\searrow 0$, apart from partially being addressed by numerical considerations
(\cite{LWZ}), seems to lack a rigorous answer up to now.\abs
The goal of the present work consists in establishing a first result in this direction,
with a main focus being on deriving an approach robust enough so as to be not necessarily restricted
to the prototypical system (\ref{KS}), but rather capable of adequately treating more complex types of interaction,
possibly also with further components.
In order to include an example for the latter which appears to be of increasing interest in the recent literature,
we shall address this problem in the context of the class of Keller-Segel systems possibly coupled to the
(Navier-)Stokes equations from fluid mechanics, and
for a fixed number $T>0$ and arbitrary $\eps>0$, we will accordingly be concerned with solutions to the class of
systems given by
\be{0eps}
	\left\{ \begin{array}{lcll}
	\,  \,  \partial_t \neps + \ueps \cdot \nabla\neps
	&=& \Delta \neps - \nabla \cdot \Big( \neps S(x,\neps,\ceps)\cdot\nabla \ceps\Big) + f(x,\neps,\ceps),
	\quad & x\in \Omega, \ t\in (0,T), \\[1mm]
	\eps \partial_t \ceps + \ueps\cdot\nabla \ceps
	&=& \Delta \ceps - \ceps + \neps,
	\quad & x\in \Omega, \ t\in (0,T), \\[1mm]
	\, \, \partial_t \ueps + \kappa (\ueps\cdot\nabla)\ueps
	&=&  \Delta \ueps + \nabla \Peps + \neps \nabla\phi,
	\qquad \nabla\cdot \ueps=0,
	\quad & x\in \Omega, \ t\in (0,T), \\[1mm]
	& & \hspace*{-41mm}
	\big(\nabla\neps - \neps S(x,\neps,\ceps)\cdot\nabla\ceps\big)\cdot \nu= \frac{\partial\ceps}{\partial\nu}=0,
	\quad \ueps=0,
	\quad & x\in\pO, \ t\in (0,T), \\[1mm]
	& & \hspace*{-40mm}
	\neps(x,0)=n_0(x), \quad \ceps(x,0)=c_0(x), \quad \ueps(x,0)=u_0(x),
	\quad & x\in\Omega.
	\end{array} \right.
\ee
Here we shall include the possibility that alternative to the choice $f\equiv 0$, the proliferation term $f$
may e.g.~represent a logistic-type source, possibly even accounting for competition with the quantity $c$
such as typically present in taxis-type models from spatial ecology where $c$ plays the role of a second species
(\cite{junping_shi}).
Moreover, our approach will be general enough so as to allow for the chemotactic interaction in (\ref{0eps})
to be described by the action of a matrix which may contain off-diagonal entries, and thus especially
be able to account for rotational flux components such as proposed in the more recent modeling literature
(\cite{xue_othmer}) but yet understood only rudimentarily from an analytical point of view (\cite{cao_lankeit},
\cite{wang_xiang_JDE2015}, \cite{wang_xiang_JDE2016}, \cite{lsxw}, \cite{win_ct_rot_SIMA}).
Correspondigly, we shall suppose that with some $K_f>0$ and some nonincreasing $f_0:[0,\infty)\to\R$ with $f_0(0)\ge 0$,
\be{f}
	\left\{ \begin{array}{l}
	f\in C^1(\bom\times [0,\infty)^2)
	\qquad \mbox{is such that $f(x,0,c)\ge 0$ \quad for all $(x,c)\in\bom\times [0,\infty)$, \qquad and that}\\[1mm]
	f_0(n) \le f(x,n,c) \le K_f \cdot (n+1)
	\quad \mbox{for all } (x,n,c)\in\bom\times [0,\infty)^2,
	\end{array} \right.
\ee
that $S=(S_{ij})_{i,j\in\{1,...,N\}}$ is such that for all $(i,j)\in \{1,...,N\}^2$,
\be{S}
	\left\{ \begin{array}{l}
	S_{ij} \in C^2(\bom\times [0,\infty)^2),
	\qquad \mbox{and that} \\[1mm]
	|S_{ij}(x,n,c)| \le K_S
	\quad \mbox{for all } (x,n,c) \in \bom\times [0,\infty)^2
	\end{array} \right.
\ee
with a positive constant $K_S$, and that apart from that the parameter $\kappa$ is any real number and the
gravitational potential in (\ref{0eps}) satisfies
\be{phi}
	\phi \in W^{2,\infty}(\Omega).
\ee
As for the initial data, our standing assumptions will be that
\be{init}
	\left\{ \begin{array}{l}
	n_0\in W^{1,\infty}(\Omega)
	\quad \mbox{is nonnegative with } n_0\not\equiv 0, \\[1mm]
	c_0\in W^{1,\infty}(\Omega)
	\quad \mbox{is nonnegative, \, and that} \\[1mm]
	u_0\in W^{2,\infty}(\Omega;\R^N)
	\quad \mbox{satisfies $\nabla\cdot u_0 \equiv 0$ and $u_0|_{\pO}=0$.}
	\end{array} \right.
\ee
\abs
Our plan is to firstly derive a general result on convergence of solutions to (\ref{0eps}) to solutions of the
associated parabolic-elliptic counterpart, and to secondly concretize this in the framework of two particular examples.
We shall thereby obtain corresponding approximation results both for certain small-data solutions to an unforced
chemotaxis-Navier-Stokes system, and for arbitrary solutions to a one-dimensional fluid-free logistic Keller-Segel model,
where as a by-product, the latter outcome will imply an apparently new result
on spontaneous emergence of arbitrarily large densities in the fully parabolic case for suitably small $\eps>0$.\abs
{\bf Main results I. A general statement on the limit $\eps\searrow 0$ in (\ref{0eps}).} \quad
Accordingly, we shall first examine the relationship between solutions to (\ref{0eps}) and those to
\be{0e}
	\left\{ \begin{array}{lcll}
	n_t + u\cdot\nabla n &=& \Delta n - \nabla \cdot (nS(x,n,c)\cdot\nabla c) + f(x,n,c),
	\qquad & x\in \Omega, \ t\in (0,T), \\[1mm]
	\qquad    u\cdot\nabla c &=& \Delta c-c+n,
	\qquad & x\in \Omega, \ t\in (0,T), \\[1mm]
	u_t + \kappa (u\cdot\nabla) u &=& \Delta u + \nabla P + n\nabla\phi,
	\qquad \nabla\cdot u=0,
	\qquad & x\in \Omega, \ t\in (0,T), \\[1mm]
	& & \hspace*{-34mm}
	\big(\nabla n - n S(x,n,c) \cdot\nabla c\big) \cdot \nu = \frac{\partial c}{\partial\nu}=0,
	\quad u=0,
	\qquad & x\in\pO, \ t\in (0,T), \\[1mm]
	& & \hspace*{-33mm}
	n(x,0)=n_0(x), u(x,0)=u_0(x),
	\qquad & x\in\Omega,
	\end{array} \right.
\ee
in a setting as general as possible. Our main result in this respect identifies a condition, yet on a given
family of solutions to (\ref{0eps}) itself, as sufficient for strong, and especially a.e.~pointwise convergence,
in the following sense.
\begin{theo}\label{theo47}
  Let $N\ge 1$ and $\Omega\subset\R^N$ be a bounded convex domain with smooth boundary, and assume that (\ref{init}) holds,
  that $\kappa\in\R$, and that $f$, $\phi$ and $S$ comply with (\ref{f}), (\ref{phi}) and  (\ref{S}).  Furthermore, suppose
  that $(\eps_j)_{j\in\N} \subset (0,\infty)$ is such that $\eps_j\searrow 0$ as $j\to\infty$, and that for some $T>0$,
  $((\neps,\ceps,\ueps,\Peps))_{\eps\in (\eps_j)_{j\in\N}}$ is such that for each  $\eps\in (\eps_j)_{j\in\N}$,
  $(\neps,\ceps,\ueps,\Peps)$ solves (\ref{0eps}) classically in $\Omega\times (0,T)$ with  $\neps\ge 0$ and $\ceps\ge 0$
  in $\Omega\times (0,T)$, and such that
  \be{47.1}
	\sup_{\eps\in (\eps_j)_{j\in\N}} \|\nabla\ceps\|_{L^\lambda((0,T);L^q(\Omega))} < \infty
  \ee
  as well as
  \be{47.2}
	\sup_{\eps\in (\eps_j)_{j\in\N}} \|\ueps\|_{L^\infty((0,T);L^r(\Omega))} < \infty
  \ee
  with some $\lambda\in (2,\infty]$, $q>N$ and $r>\max\big\{2,N\big\}$ satisfying
  \be{47.22}
	\frac{1}{\lambda} + \frac{N}{2q} < \frac{1}{2}.
  \ee
  Then there exist a subsequence $(\eps_{j_k})_{k\in\N}$ of $(\eps_j)_{j\in\N}$ and a classical solution
  $(n,c,u,P)$ of (\ref{0e}) in $\Omega \times (0,T)$ with the properties that
  \begin{eqnarray}
	& & \neps \to n
	\qquad \mbox{in } \, C^0\big(\bom\times [0,T]\big),
	\label{47.3} \\[1mm]
	& & \neps \wto n
	\qquad \mbox{in } \,  L^2\big((0,T);W^{1,2}(\Omega)\big),
	\label{47.4} \\[1mm]
	& & \ceps \to c
	\qquad \mbox{in } \,  L^\infty_{loc}\big((0,T];C^0(\bom)\big) \cap L^2_{loc}\big((0,T];W^{1,2}(\Omega)\big),
	\label{47.5} \\[1mm]
	& & \nabla\ceps \wsto \nabla c
	\qquad \mbox{in } \,  \bigcap_{\hatq>N} L^\infty\big((0,T);W^{1,\hatq}(\Omega)\big)
	\cap L^\infty\big((\Omega\times (0,T)\big)
	\qquad \mbox{and}   \label{47.6} \\[1mm]
	& & \ueps \to u
	\qquad \mbox{in }  \, C^0\big(\bom\times [0,T]; \mathbb{R}^{ N}\big)
	\cap C^{2,1}_{loc}\big(\bom\times (0,T]; \mathbb{R}^{N}\big)
	\label{47.7}
  \eea
  as $\eps=\eps_{j_k}\searrow 0$.
\end{theo}
{\bf Main results II. The fast signal diffusion limit for small-data solutions to a Keller-Segel-Navier-Stokes system.} \quad
As a first application of the latter, let us consider the case when $f\equiv 0$ in the Keller-Segel-Navier-Stokes
system (\ref{0eps}) in arbitrary spatial dimensions $N\ge 2$.
Then in light of well-known results on taxis-driven blow-up of some solutions to both the fully parabolic problem
(\ref{0eps}) as well as its parabolic-elliptic counterpart (\ref{0e}) already in the simple case $u\equiv 0$
(\cite{nagai2001}, \cite{herrero_velazquez}, \cite{win_JMPA}), regular behavior throughout the arbitrary time
interval $(0,T)$ can be expected only under appropriate additional assumptions on the initial data.
In deriving the following consequence of Theorem \ref{theo47} on this particular system, we shall accordingly
restrict our considerations to solutions emanating from suitably small initial data. In this context we will
see the following.
\begin{theo}\label{theo62}
  Let $N\ge 2$ and $\Omega\subset\R^N$ be a bounded convex domain with smooth boundary, let $\kappa\in\R$,   $p>N$, $q>N$ and $r>N$, and assume that  (\ref{phi}) and (\ref{S}) are valid.
  Then there exists $\delta=\delta(p,q,r)>0$ with the property that whenever $n_0, c_0$ and $u_0$ comply with (\ref{init})
  and satisfy
  \be{62.1}
	\|n_0\|_{L^p(\Omega)} \le \delta,
	\qquad
	\|\nabla c_0\|_{L^q(\Omega)} \le \delta
	\qquad \mbox{and} \qquad
	\|u_0\|_{L^r(\Omega)} \le \delta,
  \ee
  for all $\eps>0$ the problem (\ref{0eps}) with $f\equiv 0$ possesses a global classical solution $(\neps,\ceps,\ueps,\Peps)$.  Moreover, there exists a global classical solution $(n,c,u,P)$ of (\ref{0e}) with $f\equiv 0$ such that given any  $(\eps_j)_{j\in\N} \subset (0,\infty)$ satisfying $\eps_j\searrow 0$ as $j\to\infty$ and  each $T>0$,   one can extract
  a subsequence $(\eps_{j_k})_{k\in\N}$ fulfilling  (\ref{47.3})-(\ref{47.7}) as $\eps=\eps_{j_k}\searrow 0$.
\end{theo}
{\bf Main results III. A growth phenomenon in a fully parabolic one-dimensional Keller-Segel system with logistic source.}
\quad
As a second application of our general theory, we shall consider the family of fluid-free one-dimensional
Keller-Segel systems with logistic source, as given by
\be{L}
	\left\{ \begin{array}{ll}
	\, \,  n_{\varepsilon t} =D n_{\varepsilon xx} - (n_{\varepsilon}c_{\varepsilon x})_x + an_{\varepsilon} - bn_{\varepsilon}^2,
	\qquad & x\in (0,1), \ t>0, \\[1mm]
	\eps c_{\varepsilon t} = c_{\varepsilon xx} - c_{\varepsilon} + n_{\varepsilon},
	\qquad & x\in (0,1) , \ t>0, \\[1mm]
	\, \, n_{\varepsilon x}(0,t)=n_{\varepsilon x}(1,t)=c_{\varepsilon x}(0,t)= c_{\varepsilon x}(1,t)=0,
	\qquad &   t>0, \\[1mm]
	\, \,  n_{\varepsilon}(x,0)=n_0(x), \quad
	c_{\varepsilon}(x,0)=c_0(x),
	\qquad & x\in (0,1),
	\end{array} \right.
\ee
for $\eps>0$,  with $a\in\R$ and $b\ge 0$, and with nonnegative functions $n_0\in  W^{1,\infty}\big((0,1)\big)$ and $c_0\in W^{1,\infty}\big((0,1)\big)$.  We note that upon replacing  $n_{\varepsilon}$ by  $\tilde n_{\varepsilon}(x,\tilde t):=n_{\varepsilon}(x,t)$  with $\tilde t:=Dt$ for $(x,t)\in [0,1]\times [0,\infty)$,
this problem indeed takes the form (\ref{0eps}), and that as a well-known fact,
for each $\eps>0$ there exists a global classical solution
$(\neps,\ceps)$ for which both $\neps$ and $\ceps$ are nonnegative and bounded throughout $(0,1)\times (0,\infty)$
(\cite{win_CPDE2010}).
In fact, we shall firstly see that as $\eps\searrow 0$, these solutions approach solutions to the corresponding
parabolic-elliptic counterpart, namely, the problem
\be{Le}
	\left\{ \begin{array}{ll}
	n_t =D n_{xx} - (nc_x)_x + an - bn^2,
	\qquad & x\in (0,1), \ t>0, \\[1mm]
	\, 0  \,  = c_{xx} - c + n,
	\qquad & x\in(0,1), \ t>0, \\[1mm]
	\, n_x(0,t)=n_x(1,t)=c_x(0,t)=c_x(1,t)=0,
	\qquad &   t>0, \\[1mm]
	\, n(x,0)=n_0(x),
	\qquad & x\in (0,1),
	\end{array} \right.
\ee
in the following sense:
\begin{theo}\label{theo66}
  Let  $D>0$, $a\in\R$ and $b\ge 0$, and suppose that $n_0$ and $c_0$ are such that
  (\ref{init}) holds.
  Then for all $T>0$, the solutions $(\neps,\ceps)$ of (\ref{L}) have the property that as $\eps\searrow 0$,
  (\ref{47.3})-(\ref{47.6}) hold with the unique classical solution
  $(n,c) \in \big(C^0\big([0,1]\times [0,T]\big) \cap C^{2,1}\big([0,1]\times (0,T)\big)\big) \times C^{2,0}\big([0,1]\times (0,T)\big)$
  of (\ref{Le}).
\end{theo}
Building on this result, we shall secondly discover that solutions to the fully parabolic problem (\ref{L})
can spontaneously generate arbitrarily large densities, possibly at intermediate time scales, provided that the
parameters $D$ and $\eps$ satisfy appropriate smallness conditions:
\begin{theo}\label{theo67}
  Let $a\in\R$ and $b\in [0,1)$. Then there exist $T>0$ and a nonnegative function $n_0\in W^{1,\infty}\big((0,1)\big)$ with the following property: For all $M>0$ one can find $D_0>0$ such that for each $D\in (0,D_0)$ and any nonnegative $c_0\in W^{1,\infty}\big((0,1)\big)$  there exist $x_0\in (0,1), t_0\in (0,T)$ and $\eps_0>0$ such that for any choice of $\eps\in (0,\eps_0)$, the corresponding
  solution $(\neps,\ceps)$ of (\ref{L}) satisfies
  \be{67.1}
	\neps(x_0,t_0) \ge M.
  \ee
\end{theo}
{\bf Key steps in our analysis.}   \quad
The crucial role of the assumptions from Theorem \ref{theo47}, and especially of the inequality (\ref{47.22}) therein,
will already become clear in Section \ref{sect2}, in which
we will derive some $\eps$-independent estimates for general solutions to ({\ref{0eps}}) under
presupposed bounds on $\nabla \ceps$ and $\ueps$ of the considered form.
Complementing these estimates by further compactness properties will allow for passing to the limit along subsequences,
with regard to the components $\neps$ and $\ueps$ already in the flavor claimed in Theorem \ref{theo47};
as for $\ceps$, however, due to lacking uniform parabolicity in the equation describing its evolution we will at that
stage only be able to conclude a weak convergence property in $L^2((0,T);W^{1,2}(\Omega))$.\abs
A key step will thereafter consist in improving this knowledge, which will be achieved through several steps:
After firstly showing that the limit $c$ satisfies its respective subproblem of (\ref{0e}) in a weak sense,
we can exploit the correspondingly satisfies integral identity to successively
establish H\"older regularity of $c$, $\nabla c$ and $D^2 c$ in Section \ref{sect5.1}.
The main step will then be accomplished by
ensuring $L^2$ integrability of $c_t$, locally away from the temporal origin, in Section \ref{sect5.2}.
Our derivation thereof will rely on suitably estimating the difference quotients
\bas
	z_h(x,t):=\frac{c(x,t+h)-c(x,t)}{h}=:z^1_h(x,t)+z^2_h(x,t),
	\qquad x\in\Omega, \ t\in (\tau,T-h_0),
  \eas
where $z^1_h$ and $z^2_h$ denote the classical solution of two linear elliptic equations, whose forcing terms
involve the time derivatives of $n$ and $u$.
Thanks to the availability of appropriate regularity information on the latter, by utilizing
standard elliptic regularity theory we will infer that indeed
$c_t$ belongs to $L^2_{loc}(\bom\times (0,T])$. This in turn will allow us to adequately control the difference
$\ceps-c$ through analyzing a parabolic equation therefor, and hence verify Theorem \ref{theo47} in Section \ref{sect6}.\abs
Sections \ref{sect_t62} and \ref{sect8} will thereafter be devoted to the proofs of Theorem \ref{theo62} and of
Theorems \ref{theo66} and \ref{theo67}, respectively.
\mysection{Some general estimates}\label{sect2}
In this section we collect some estimates which are valid for general solutions to systems of the form (\ref{0eps}),
and which are independent of the particular choice of $\eps>0$, partially under presupposed bounds resembling
those in (\ref{47.1}) and (\ref{47.2}).
These estimates will firstly be used as a fundament for our proof of Theorem \ref{theo47}, and secondly some of them
will afterwards serve as helpful ingredients for the derivation of Theorem \ref{theo62} in Section \ref{sect_t62}.\abs
Let us start with a fairly evident observation.
\begin{lem}\label{lem31}
  Suppose that   (\ref{f}), (\ref{phi}) and (\ref{init}) hold, and let $\kappa\in\R$ and $T>0$.
  Then there exists $C=C(T)>0$ such that whenever $(\neps, \ceps, \ueps, \Peps)$
  is a classical solution of (\ref{0eps}) in $\Omega\times (0,T)$  for some $\eps>0$, we have $\neps\ge 0$ and $\ceps\ge 0$ in $\Omega\times (0,T)$ as well as
  \be{mass}
	\|\neps(\cdot,t)\|_{L^1(\Omega)} \le C
	\qquad \mbox{for all } t\in (0,T).
  \ee
\end{lem}
\proof
  According to the lower bound for $f_0$ and hence for $f$ in (\ref{f}), nonnegativity of $\neps$ results from an application  of the maximum principle to the first equation in (\ref{0eps}). In view of the second equation therein, by the same token  this in turn entails nonnegativity also of $\ceps$.\abs
  Next, integrating the first equation in (\ref{0eps}) shows that since $\nabla\cdot \ueps\equiv 0$, due to the upper bound for
  $f$ from (\ref{f}) we have
  \bas
	\frac{d}{dt} \io \neps = \io f(x,\neps, \ceps) \le K_f \io \neps + K_f |\Omega|
	\qquad \mbox{for all } t\in (0,T),
  \eas
  from which (\ref{mass}) readily results upon a time integration.
\qed
The next lemma already makes full use of supposedly present bounds in the style of (\ref{47.1}) and (\ref{47.2}),
and especially of the relation (\ref{47.22}) involving the parameters therein.
\begin{lem}\label{lem32}
 Suppose that   (\ref{f}),  (\ref{phi}), (\ref{S}) and  (\ref{init}) hold, and let $\kappa\in\R$.
  Then for all $T>0$,  $L>0$, $\lambda\in (2,\infty], q>N, r>N$ such that (\ref{47.22}) holds, there exists  $C=C(T,\lambda,q,r,K_S,L)>0$ such that whenever  $(\neps, \ceps, \ueps, \Peps)$ is a classical solution of (\ref{0eps})  in $\Omega\times (0,T)$  for some $\eps>0$  fulfilling
  \be{32.01}
	\|\nabla \ceps\|_{L^\lambda((0,T);L^q(\Omega))} \le L
  \ee
  and
  \be{32.02}
	\|\ueps(\cdot,t)\|_{L^r(\Omega)} \le L
	\qquad \mbox{for all } t\in (0,T),
  \ee
  we have
  \be{32.1}
	\|\neps(\cdot,t)\|_{L^\infty(\Omega)} \le C
	\qquad \mbox{for all } t\in (0,T).
  \ee
\end{lem}
\proof
Omitting the subscript $\eps$ for notational convenience,  without loss of generality assuming that $\lambda<\infty$ and following an essentially well-established procedure  (cf.~e.g.~\cite{taowin_ZAMP}), we  estimate
  \bas
	M(T'):=\sup_{t\in (0,T')} \|n(\cdot,t)\|_{L^\infty(\Omega)},
	\qquad T'\in (0,T)
  \eas
  by representing $n$ via an associated Duhamel formula.   Indeed, using the maximum principle and (\ref{f}) as well as known smoothing properties of the Neumann heat semigroup
  $(e^{t\Delta})_{t\ge 0}$ in $\Omega$ we see that fixing any $\mu_1=\mu_1(\lambda,q) \in (N,q)$ and $\mu_2=\mu_2(r)\in (N,r)$
  such that $\frac{1}{\lambda}+\frac{N}{2\mu_1}<\frac{1}{2}$, with some $C_1=C_1(\lambda,q,r)>0$ we have
  \bea{32.2}
	\hspace*{-8mm}
	n(\cdot,t)
	&=& e^{t\Delta} n_0
	- \int_0^t e^{(t-s)\Delta} \nabla \cdot \Big( n(\cdot,s) S(\cdot,n(\cdot,s),c(\cdot,s)) \cdot\nabla c(\cdot,s)\Big) ds
	\nn\\
	& & - \int_0^t e^{(t-s)\Delta} \nabla \cdot (n(\cdot,s) u(\cdot,s)) ds
	+ \int_0^t e^{(t-s)\Delta} f(\cdot,n(\cdot,s), c(\cdot,s)) ds \nn\\[1mm]
	&\le& \|n_0\|_{L^\infty(\Omega)}
	+ C_1 \int_0^t (t-s)^{-\frac{1}{2}-\frac{N}{2\mu_1}}
		\|n(\cdot,s) S(\cdot,n(\cdot,s),c(\cdot,s)) \cdot \nabla c(\cdot,s)\|_{L^{\mu_1}(\Omega)} ds \nn\\
	& & + C_1 \int_0^t (t-s)^{-\frac{1}{2}-\frac{N}{2\mu_2}} \|n(\cdot,s)u(\cdot,s)\|_{L^{\mu_2}(\Omega)} ds
	+ C_1 \int_0^t (t-s)^{-\frac{1}{2}} \|n(\cdot,s)+1\|_{L^N(\Omega)} ds
  \eea
  for all $t\in (0,T)$.
  Here using (\ref{S}) and the H\"older inequality along with (\ref{mass}) and our hypotheses (\ref{32.01}) and (\ref{32.02}),
  we find positive constants $C_2=C_2(K_S), C_3=C_3(K_S), C_4=C_4(T, K_S),  C_5=C_5(T, r,L)$ and $C_6=C_6(T)$ such that
  \bas
	\|n(\cdot,s) S(\cdot,n(\cdot,s),c(\cdot,s))\cdot\nabla c(\cdot,s)\|_{L^{\mu_1}(\Omega)}
	&\le& C_2 \|n(\cdot,s)\|_{L^\frac{q\mu_1}{q-\mu_1}(\Omega)} \|\nabla c(\cdot,s)\|_{L^q(\Omega)} \\
	&\le& C_3 \|n(\cdot,s)\|_{L^\infty(\Omega)}^{a_1} \|n(\cdot,s)\|_{L^1(\Omega)}^{1-a_1}
	\|\nabla c(\cdot,s)\|_{L^q(\Omega)} \\
	&\le& C_4 M^{a_1}(T') \|\nabla c(\cdot,s)\|_{L^q(\Omega)}
	\qquad \mbox{for all } s\in (0,T')
  \eas
  and
  \bas
	\|n(\cdot,s) u(\cdot,s)\|_{L^{\mu_2}(\Omega)}
	&\le& \|n(\cdot,s)\|_{L^\frac{r\mu_2}{r-\mu_2}(\Omega)} \|u(\cdot,s)\|_{L^r(\Omega)} \\
	&\le& \|n(\cdot,s)\|_{L^\infty(\Omega)}^{a_2} \|n(\cdot,s)\|_{L^1(\Omega)}^{1-a_2} \|u(\cdot,s)\|_{L^r(\Omega)} \\
	&\le& C_5 M^{a_2}(T')
	\qquad \mbox{for all } s\in (0,T')
  \eas
  as well as
  \bas
	\|n(\cdot,s)+1\|_{L^N(\Omega)}
	&\le& \|n(\cdot,s)+1\|_{L^\infty(\Omega)}^{a_3} \|n(\cdot,s) +1 \|_{L^1(\Omega)}^{1-a_3} \\
	&\le& C_6 M^{a_3}(T') + C_6
	\qquad \mbox{for all } s\in (0,T')
  \eas
  with $a_1:=\frac{q\mu_1-q+\mu_1}{q\mu_1}\in (0,1)$, $a_2:=\frac{ru_2-r+\mu_2}{r\mu_2}\in (0,1)$
  and $a_3:=\frac{N-1}{N}\in (0,1)$.
  Since the inequalities $\frac{1}{\lambda} + \frac{N}{2\mu_1}<\frac{1}{2}$ and $\mu_2>N$ moreover warrant that
  $(\frac{1}{2}+\frac{N}{2\mu_1})\cdot \frac{\lambda}{\lambda-1}<1$ and $\frac{1}{2} + \frac{N}{2\mu_2}<1$, by using the
  H\"older inequality two more times we thus infer from (\ref{32.2}) and the nonnegativity of $n$ that there exist
  $C_8=C_8(T,\lambda,q,r,K_S,L)>0$, $C_9=C_9(T, \lambda,q)>0$ and $C_{10}=C_{10}(r)>0$ such that
  \bas
	\|n(\cdot,t)\|_{L^\infty(\Omega)}
	&\le& C_8 + C_8 M^{a_1}(T') \cdot \bigg\{
	\int_0^t (t-s)^{-(\frac{1}{2}+\frac{N}{2\mu_1}) \cdot \frac{\lambda}{\lambda-1}} ds \bigg\}^\frac{\lambda-1}{\lambda}
	\cdot \bigg\{ \int_0^t \|\nabla c(\cdot,s)\|_{L^q(\Omega)}^\lambda ds \bigg\}^\frac{1}{\lambda} \nn\\
	& & + C_8 M^{a_2}(T') \cdot \int_0^t (t-s)^{-\frac{1}{2} - \frac{N}{2\mu_2}} ds   + C_8 M^{a_3}(T') \nn\\[1mm]
	&\le& C_8 + C_8 C_9 M^{a_1}(T')
	+ C_8 C_{10} M^{a_2}(T') + C_8 M^{a_3}(T')
	\qquad \mbox{for all } t\in (0,T').
  \eas
  Hence, by Young's inequality,
  \bas
	M(T') \le C_{11} + C_{11} M^a(T')
	\qquad \mbox{for all } T'\in (0,T),
  \eas
  where $a:=\max\{a_1,a_2,a_3\}$ satisfies $a\in (0,1)$, and where $C_{11}:=2C_8+C_8 C_9 + C_8 C_{10}$.
  As therefore
  \bas
	M(T') \le \max \Big\{ 1 \, , \, (2C_{11})^\frac{1}{1-a} \Big\}
	\qquad \mbox{for all } T'\in (0,T),
  \eas
  we have thus established (\ref{32.1}).
\qed
As a consequence of the latter estimate for $\neps$, by means of quite a similar argument, essentially well-established
in the theory of the Navier-Stokes system, we can again use the boundedness assumption (\ref{32.02}) in order to appropriately
control the fluid velocity field as follows.
\begin{lem}\label{lem33}
   Suppose that   (\ref{f}),  (\ref{phi}), (\ref{S}) and  (\ref{init}) hold, and let $\kappa\in\R$.
  Then for all $T>0$,  $L>0$, $\lambda\in (2,\infty], q>N, r>\max\{2,N\}$  fulfilling (\ref{47.22}), there exist
  $\alpha=\alpha(r)\in (\frac{1}{2},1)$, $\rho=\rho(r)> \max\big\{1,\frac{N}{2\alpha} \big\}$, $\theta=\theta(r)\in (0,1)$ and   $C=C(T,\lambda,q,r,K_S,\kappa,L)>0$ such that  if  for some  $\eps>0$  $(\neps, \ceps, \ueps, \Peps)$ is a classical solution of (\ref{0eps}) in $\Omega\times (0,T)$ for which (\ref{32.01}) and (\ref{32.02}) are valid, then
  \be{33.1}
	\|A^\alpha \ueps(\cdot,t)\|_{L^\rho(\Omega)} \le C
	\qquad \mbox{for all } t\in (0,T)
  \ee
  and
  \be{33.2}
	\|\ueps\|_{C^{\theta,\frac{\theta}{2}}(\bom\times [0,T])} \le C.
  \ee
\end{lem}
\proof
  Since $r>N$ and thus $\frac{1}{2} + \frac{N}{2r}<1$, it is possible to fix $\alpha=\alpha(r)\in (\frac{1}{2},1)$ close to
  $\frac{1}{2}$ such that
  \be{33.22}
	\alpha+\frac{N}{2r}<1,
  \ee
  and thereafter take $\beta\in (\frac{1}{2},1)$ such that $\beta<\alpha$.
  Then using that $\frac{N}{2\alpha}<N<r$ and that also $\frac{r}{r-1}<r$, we can pick
  $\rho=\rho(r)>\max \big\{1,\frac{N}{2\alpha}\big\}$ such that $\rho \le r$ and $\rho>\frac{r}{r-1}$, observing that the latter  ensures that $\mu=\mu(r):=\frac{r\rho}{r+\rho}$ satisfies $\mu>1$.
 Again dropping the index $\eps$ and  as moreover $\mu<\rho$, relying on a variation-of-constants representation of $u$ we may  employ known smoothing properties  of the Stokes semigroup $(e^{-tA})_{t\ge 0}$ (\cite{giga1986}) to find $C_1=C_1(r,\kappa)>0$ such that
  \bea{33.3}
	\|A^\alpha u(\cdot,t)\|_{L^\rho(\Omega)}
	&=& \Bigg\| A^\alpha e^{-tA} u_0
	- \kappa \int_0^t A^\alpha e^{-(t-s)A} \proj \Big[ (u(\cdot,s)\cdot\nabla) u(\cdot,s) \Big] ds \nn\\
	& & \hspace*{20mm}
	+ \int_0^t A^\alpha e^{-(t-s)A} \proj [n(\cdot,s)\nabla \phi] ds \Bigg\|_{L^\rho(\Omega)} \nn\\
	&\le& \|A^\alpha u_0\|_{L^\rho(\Omega)}
	+ C_1 \int_0^t (t-s)^{-\alpha-\frac{N}{2}(\frac{1}{\mu}-\frac{1}{\rho})}
	\|(u(\cdot,s)\cdot\nabla) u(\cdot,s)\|_{L^\mu(\Omega)} ds \nn\\
	& & + C_1 \int_0^t (t-s)^{-\alpha} \|n(\cdot,s)\|_{L^\rho(\Omega)} ds
	\qquad \mbox{for all } t\in (0,T).
  \eea
  Here since clearly $\mu<r$, we can employ the H\"older inequality to see that thanks to (\ref{32.02}) and the inequalities
  $\alpha>\beta>\frac{1}{2}$ and $\rho\le r$, the continuity of the embedding $D(A_\rho^\beta) \hra W^{1,\rho}(\Omega;\R^N)$
  {(\cite{friedman}  \cite{henry}}) and a well-known interpolation property guarantee that with some
  $C_2=C_2(r)>0, C_3=C_3(r)>0$ and $C_4=C_4(r)>0$ we have
  \bas
	\|(u(\cdot,s)\cdot\nabla )u(\cdot,s)\|_{L^\mu(\Omega)}
	&\le& \|u(\cdot,s)\|_{L^r(\Omega)} \|\nabla u(\cdot,s)\|_{L^\frac{r\mu}{r-\mu}(\Omega)} \\
	&\le& L \|\nabla u(\cdot,s)\|_{L^\rho(\Omega)} \\
	&\le& C_2 L \|A^\beta u(\cdot,s)\|_{L^\rho(\Omega)} \\
	&\le& C_3 L \|A^\alpha u(\cdot,s)\|_{L^\rho(\Omega)}^a \|u(\cdot,s)\|_{L^\rho(\Omega)}^{1-a} \\
	&\le& C_4 L \|A^\alpha u(\cdot,s)\|_{L^\rho(\Omega)}^a \|u(\cdot,s)\|_{L^r(\Omega)}^{1-a} \\
	&\le& C_4 L^{2-a} M^a(T')
	\qquad \mbox{for all $s\in (0,T')$ and any $T'\in (0,T)$}
  \eas
  if we let $a:=\frac{\beta}{\alpha}\in (0,1)$ and
  \bas
	M(T'):=\sup_{t\in (0,T')} \|A^\alpha u(\cdot,t)\|_{L^\rho(\Omega)},
	\qquad T'\in (0,T).
  \eas
  As Lemma \ref{lem32} in particular implies the existence of $C_5=C_5(T, \lambda, q, r, K_S, L)>0$ such that
  \bas
	\|n(\cdot,t)\|_{L^\rho(\Omega)} \le C_5
	\qquad \mbox{for all } t\in (0,T),
  \eas
  noting that $\alpha<1$ and that
  \bas
	\alpha+ \frac{N}{2}\Big(\frac{1}{\mu}-\frac{1}{\rho}\Big)
	= \alpha+\frac{N}{2}\Big( \frac{r+\rho}{r\rho} - \frac{1}{\rho}\Big)
	= \alpha+\frac{N}{2r}<1
  \eas
  by (\ref{33.22}), we thus conclude from (\ref{33.3}) and (\ref{init}) that there exists $C_6=C_6(T, \lambda, q, r, K_S, \kappa, L)>0$  such that
  \bas
	M(T') \le C_6 + C_6 M^a(T')
	\qquad \mbox{for all $T'\in (0,T)$},
  \eas
  which implies (\ref{33.1}) due to the fact that $a<1$.\abs
  Now by a straightforward adaptation of a well-known reasoning (\cite{friedman}), in quite a similar manner it is
  furthermore possible to find $\theta_1=\theta_1(r)\in (0,1)$ and $C_7=C_7(T, \lambda, q, r, K_S, \kappa, L)>0$ fulfilling
  \bas
	\|A^\alpha u(\cdot,t) - A^\alpha u(\cdot,t_0)\|_{L^\rho(\Omega)} \le C_7 |t-t_0|^{\theta_1}
	\qquad \mbox{for all $t\in (0,T)$ and } t_0\in (0,T),
  \eas
  which finally implies (\ref{33.2}) due to the fact that $D(A_\rho^\alpha) \hra C^{\theta_2}(\bom;\R^N)$
  for any $\theta_2$ from the nonempty interval $\big(0,2\alpha-\frac{N}{\rho}\big)$ (\cite{henry}).
\qed
Let us finally prepare an argument that will, before becoming substantial for the derivation of Theorem \ref{theo62} in
Section \ref{sect_t62}, inter alia reveal in Lemma \ref{lem50} that
the assumptions (\ref{47.1}) and (\ref{47.2}) actually imply boundedness of $\nabla c$ in
$L^\infty((0,T);L^{\hatq}(\Omega))$ for arbitrarily large $\hatq$.
The following lemma is the only place in this paper where convexity of $\Omega$ is explicitly needed.
\begin{lem}\label{lem01}
%
For all $p>\max\{N,2\}$, $q\ge 2$ and $r\in (2,\infty]$ there exists $C=C(p,q)>0$ such that
  if $(\neps, \ceps, \ueps, \Peps)$ is a classical solution of (\ref{0eps}) in $\Omega\times (0,T)$  for some $\eps>0$ and $T>0$, then
  \bea{01.1}
	& & \hspace*{-20mm}
	\frac{\eps}{q} \frac{d}{dt} \io |\nabla \ceps|^q
	+ \frac{1}{4} \io |\nabla \ceps|^{q-2} |D^2 \ceps|^2
	+ \Big(1-\frac{1}{q^2}\Big) \io |\nabla \ceps|^q \nn\\
	&\le& C \|\neps\|_{L^p(\Omega)}^q
	+ C\|\ueps\|_{L^r(\Omega)}^2 \|\nabla \ceps\|_{L^\frac{qr}{r-2}(\Omega)}^q
	\qquad \mbox{for all } t\in (0,T),
  \eea
  where we interpret $\frac{qr}{r-2}$ as coinciding with $q$ if $r=\infty$.
\end{lem}
\proof
Once more omitting the subscript $\eps$ for convenience,   by means of the second equation in (\ref{0eps}), we see that  for all $t\in (0,T)$,
  \bea{01.2}
	\frac{\eps}{q} \frac{d}{dt} \io |\nabla c|^q
	&=& \io |\nabla c|^{q-2} \nabla c \cdot \nabla \Big\{ \Delta c - c + n - u\cdot\nabla c \Big\} \nn\\
	&=& \frac{1}{2} \io |\nabla c|^{q-2} \Delta |\nabla c|^2
	- \io |\nabla c|^{q-2} |D^2 c|^2 \nn\\
	& & - \io |\nabla c|^q
	+ \io |\nabla c|^{q-2} \nabla c\cdot\nabla n
	- \io |\nabla c|^{q-2} \nabla c \cdot \nabla (u\cdot\nabla c) \nn\\
	&\le& - \io |\nabla c|^{q-2} |D^2 c|^2
	- \io |\nabla c|^q \nn\\
	& \quad & - \io n|\nabla c|^{q-2} \Delta c
	-  (q-2)  \io n|\nabla c|^{q-4}  \nabla c \cdot  (D^2 c\cdot\nabla c) \nn\\
	& & + \io (u\cdot\nabla c) |\nabla c|^{q-2} \Delta c
	+ (q-2) \io (u\cdot\nabla c) |\nabla c|^{q-4} \nabla c \cdot (D^2 c\cdot\nabla c),
  \eea
  because  of  $\frac{\partial |\nabla c|^2}{\partial\nu}\le 0$ on $\pO\times (0,T)$ due to the convexity of $\Omega$  (\cite{lions_ARMA}), and because of $q\ge 2$.   Here two applications of Young's inequality and the H\"older inequality show that abbreviating $C_1(q):=\sqrt{2}+q-2$
  we have
  \bea{01.3}
	& & \hspace*{-30mm}
	- \io n|\nabla c|^{q-2} \Delta c
	-   (q-2) \io n|\nabla c|^{q-4}\nabla c \cdot (D^2 c\cdot\nabla c)    \nn\\
	&\le& C_1(q) \io n|\nabla c|^{q-2} |D^2 c| \nn\\
	&\le& \frac{1}{4} \io |\nabla c|^{q-2} |D^2 c|^2
	+ C_1^2(q) \io n^2 |\nabla c|^{q-2} \nn\\
	&\le& \frac{1}{4} \io |\nabla c|^{q-2} |D^2 c|^2
	+ C_1^2(q) \|n\|_{L^p(\Omega)}^2 \|\nabla c\|_{L^\frac{p(q-2)}{p-2}(\Omega)}^{q-2}
  \eea
  and
  \bea{01.4}
	& & \hspace*{-30mm}
	\io (u\cdot\nabla c) |\nabla c|^{q-2} \Delta c
	+ (q-2) \io (u\cdot\nabla c) |\nabla c|^{q-4} \nabla c \cdot (D^2 c\cdot\nabla c) \nn\\
	&\le& C_1(q) \io |u| |\nabla c|^{q-1} |D^2 c| \nn\\
	&\le& \frac{1}{4} \io |\nabla c|^{q-2} |D^2 c|^2
	+ C_1^2(q) \io |u|^2 |\nabla c|^q \nn\\
	&\le& \frac{1}{4} \io |\nabla c|^{q-2} |D^2 c|^2
	+ C_1^2(q) \|u\|_{L^r(\Omega)}^2 \|\nabla c\|_{L^\frac{qr}{r-2}(\Omega)}^q
  \eea
  for all $t\in (0,T)$.
  Now since $p>\max\{N,2\}$ and $q\ge 2$ ensure that $p(q-2)(N-2) < N(p-2)q$, it follows that $W^{1,2}(\Omega)$ is
  continuously embedded into $L^\frac{2p(q-2)}{(p-2)q}(\Omega)$.
  Again using Young's inequality, we can therefore find $C_2(p,q)>0$ and $C_3(p,q)>0$ such that
  \bas
	C_1^2(q) \|n\|_{L^p(\Omega)}^2 \|\nabla c\|_{L^\frac{p(q-2)}{p-2}(\Omega)}^{q-2}
	&\le& C_2(p,q) \|n\|_{L^p(\Omega)}^2 \cdot \bigg\{
	\Big\| \nabla |\nabla c|^\frac{q}{2} \Big\|_{L^2(\Omega)}^2
	+ \Big\| |\nabla c|^\frac{q}{2} \Big\|_{L^2(\Omega)}^2 \bigg\}^\frac{q-2}{q} \\
	&\le& \frac{1}{q^2} \cdot \bigg\{
	\Big\| \nabla |\nabla c|^\frac{q}{2} \Big\|_{L^2(\Omega)}^2
	+ \Big\| |\nabla c|^\frac{q}{2} \Big\|_{L^2(\Omega)}^2 \bigg\}
	+ C_3(p,q) \|n\|_{L^p(\Omega)}^q \\
	&=& \frac{1}{4} \io |\nabla c|^{q-4} |D^2 c \cdot \nabla c|^2
	+ \frac{1}{q^2} \io |\nabla c|^q
	+ C_3(p,q) \|n\|_{L^p(\Omega)}^q \\
	&\le& \frac{1}{4} \io |\nabla c|^{q-2} |D^2 c|^2
	+ \frac{1}{q^2} \io |\nabla c|^q
	+ C_3(p,q) \|n\|_{L^p(\Omega)}^q
  \eas
  for all $t\in (0,T)$.
  Therefore, (\ref{01.3}) and (\ref{01.4}) when inserted into (\ref{01.2}) show that
  \bas
	\frac{\eps}{q} \frac{d}{dt} \io |\nabla c|^q
	&+& \frac{1}{4} \io |\nabla c|^{q-2} |D^2 c|^2
	+ \Big(1-\frac{1}{q^2}\Big) \io |\nabla c|^q \\
	&\le& C_3(p,q) \|n \|_{L^p(\Omega)}^q
	+ C_1^2(q) \|u\|_{L^r(\Omega)}^2 \|\nabla c\|_{L^\frac{qr}{r-2}(\Omega)}^q
	\qquad \mbox{for all } t\in (0,T),
  \eas
  which directly results in (\ref{01.1}).
\qed
\mysection{Regularity and compactness properties implied by the hypotheses from Theorem \ref{theo47}}
Next concentrating on the particular setup created by Theorem \ref{theo47}, in this part we will augment
the estimates from the previous section by further compactness properties which will allow for passing
to the limit, already partially in the flavor claimed in Theorem \ref{theo47}.\abs
Firstly, the $L^\infty$ bound from Lemma \ref{lem32} can quite immediately be improved into an estimate
in some H\"older space by means of standard parabolic theory.
\begin{lem}\label{lem34}
  Suppose that the assumptions of Theorem \ref{theo47} are satisfied.
  Then there exist $\theta\in (0,1)$ and $C>0$ such that
  \be{34.1}
	\|\neps\|_{C^{\theta,\frac{\theta}{2}}(\bom\times [0,T])} \le C
	\qquad \mbox{for all } \eps\in (\eps_j)_{j\in\N}.
  \ee
\end{lem}
\proof
  We rewrite the first equation in (\ref{0eps}) in the form
  \bas
	\partial_t \neps = \nabla \cdot a_\eps(x,t,\nabla\neps) + b_\eps(x,t),
	\qquad x\in\Omega, \ t\in (0,T),
  \eas
  with
  \bas
	a_\eps(x,t,\xi):=\xi - \neps(x,t) S(x,\neps(x,t),\ceps(x,t)) \cdot\nabla\ceps(x,t) - \neps(x,t)\ueps(x,t),
	\ (x,t,\xi)\in \Omega\times (0,T)\times \R^N,
  \eas
  and
  \bas
	b_\eps(x,t):=f(x,\neps(x,t),\ceps(x,t)),
	\qquad (x,t)\in\Omega\times (0,T),
  \eas
  Then due to Young's inequality and (\ref{f}), Lemma \ref{lem32} and Lemma \ref{lem33} yield
  positive constants $C_1$ and $C_2$ such that for all $\eps\in (\eps_j)_{j\in\N}$,
  \bas
	a_\eps(x,t,\xi)\cdot\xi
	\ge \frac{|\xi|^2}{2} - C_1 |\nabla \ceps(x,t)|^2 - C_1
	\qquad \mbox{for all } (x,t,\xi)\in \Omega\times (0,T)\times \R^N
  \eas
  and
  \bas
	|a_\eps(x,t,\xi)| \le |\xi| + C_2 |\nabla\ceps(x,t)| + C_2
	\qquad \mbox{for all } (x,t,\xi)\in \Omega\times (0,T)\times \R^N
  \eas
  as well as
  \bas
	|b(x,t)| \le C_3
	\qquad \mbox{for all } (x,t)\in \Omega\times (0,T).
  \eas
  Since (\ref{32.01}) provides a bound for $|\nabla\ceps|^2$ in $L^\frac{\lambda}{2}((0,T);L^\frac{q}{2}(\Omega))$,
  with the exponents therein satisfying $\frac{2}{\lambda} + \frac{N}{2\cdot\frac{q}{2}}=\frac{2}{\lambda}+\frac{N}{q}<1$
  by (\ref{47.22}), the estimate (\ref{34.1}) directly results on applying a standard result on H\"older regularity
  in scalar parabolic equations (\cite[Theorem 1.3, Remark 1.4]{porzio_vespri}).
\qed
Thanks to standard Schauder estimates for the Stokes system, the latter directly entails bounds for $\ueps$
even in higher-order H\"older spaces, at least locally away from the initial time.
\begin{lem}\label{lem37}
  Under the assumptions of Theorem \ref{theo47}, for each $\tau\in (0,T)$ one can find $\theta\in (0,1)$ and $C>0$ such that
  \be{37.1}
	\|\ueps\|_{C^{2+\theta,1+\frac{\theta}{2}}(\bom\times [\tau,T])} \le C
	\qquad \mbox{for all } \eps\in (\eps_j)_{j\in\N}.
  \ee
\end{lem}
\proof
  Thanks to the estimates provided by Lemma \ref{lem32} and Lemma \ref{lem33}, this follows upon a straightforward
  application of well-known Schauder theory for the linear inhomogeneous Stokes evolution equation
  (\cite{solonnikov}).
\qed
As a further implication of Lemma \ref{lem32}, by way of a standard testing procedure we can also obtain further
bounds for the second solution component which, if the parameter $s$ therein is chosen large, may partially
go beyond the information invested through (\ref{47.1}).
\begin{lem}\label{lem48}
  Suppose that the assumptions of Theorem \ref{theo47} are satisfied with some
  $(\eps_j)_{j\in\N} \subset (0,\infty)$.		
  Then for all $s\ge 2$ there exists $C(s)>0$ such that for all $\eps\in (\eps_j)_{j\in\N}$,
  \be{48.1}
	\io \ceps^s(\cdot,t) \le C(s)
	\qquad \mbox{for all $t\in (0,T)$},
  \ee
  and that
  \be{48.2}
	\int_0^T \io \ceps^{s-2} |\nabla \ceps|^2 \le C(s).
  \ee
\end{lem}
\proof
  We multiply the second equation in (\ref{0eps}) by $\ceps^{s-1}$ and use Young's inequality to see that since
  $\nabla\cdot\ueps\equiv 0$,
  \bas
	\frac{\eps}{s} \frac{d}{dt} \io \ceps^s
	+ (s-1)\io \ceps^{s-2} |\nabla \ceps|^2
	+ \io \ceps^s
	= \io \neps \ceps^{s-1}
	\le \frac{s-1}{s} \io \ceps^s + \frac{1}{s} \io \neps^s
	\qquad \mbox{for all } t\in (0,T),
  \eas
  so that Lemma \ref{lem32} implies the existence of $C_1>0$ such that for all $\eps \in (\eps_j)_{j\in\N}$,
  $\yeps(t):=\io \ceps^s(\cdot,t), \ t\in [0,T]$, satisfies
  \bas
	\eps \yeps'(t) + \yeps(t) + s(s-1) \io \ceps^{s-2} |\nabla\ceps|^2 \le C_1
	\qquad \mbox{for all } t\in (0,T).
  \eas
  This firstly entails by a comparison argument that
  \bas
	\yeps(t) \le \max \bigg\{ \io c_0^s \, , \, C_1\bigg\}	
	\qquad \mbox{for all $t\in (0,T)$ and } \eps \in (\eps_j)_{j\in\N},
  \eas
  and secondly ensures upon integration that
  \bas
	s(s-1) \int_0^T \io \ceps^{s-2} |\nabla\ceps|^2 \le \eps \io c_0^s + C_1 T
	\qquad \mbox{for all } \eps \in (\eps_j)_{j\in\N},
  \eas
  so that both (\ref{48.1}) and (\ref{48.2}) directly follow, because $(\eps_j)_{j\in\N}$ is bounded.
\qed
Similarly, Lemma \ref{lem32} together with the latter entails an estimate for $\nabla\neps$:
\begin{lem}\label{lem10}
  Suppose that the assumptions of Theorem \ref{theo47} are satisfied with some
  $(\eps_j)_{j\in\N} \subset (0,\infty)$.		
  Then there exists $C>0$ such that
  \be{10.1}
	\int_0^T \io |\nabla\neps|^2 \le C
	\qquad \mbox{for all } \eps \in (\eps_j)_{j\in\N}.
  \ee
\end{lem}
\proof
  Using $\neps$ as a test function in the first equation in (\ref{0eps}) and relying on (\ref{S}), (\ref{f})
  as well as Young's inequality, we find that again since $\nabla\cdot\ueps\equiv 0$,
  \bas
	\frac{1}{2} \frac{d}{dt} \io \neps^2 + \io |\nabla\neps|^2
	&=& \io \neps (S(x,\neps,\ceps) \cdot \nabla\ceps) \cdot\nabla\neps
	+ \io \neps f(x,\neps,\ceps) \\
	&\le& \frac{1}{2} \io |\nabla\neps|^2
	+ \frac{K_S^2}{2} \|\neps\|_{L^\infty(\Omega)}^2 \io |\nabla\ceps|^2
	+ K_f \io \neps + K_f \io \neps^2
  \eas
  for all $t\in (0,T)$.
  In view of the bounds provided by Lemma \ref{lem32} and Lemma \ref{lem48}, upon a time integration this readily yields
  (\ref{10.1}).
\qed
To prepare a useful ingredient for our subsequent analysis concerning the time regularity of the limit $c$ to be obtained,
we note the following weak but eventually helpful regularity information on $\partial_t \neps$.
For its formulation and for later reference, let us agree on using the abbreviation
$W^{2,2}_N(\Omega):=\{\psi\in W^{2,2}(\Omega) \ | \ \frac{\partial\psi}{\partial\nu}=0 \mbox{ on } \pO\}$.
\begin{lem}\label{lem36}
  Under the assumptions from Theorem \ref{theo47}, there exists $C>0$ such that
  \be{36.1}
	\int_0^T \|\partial_t \neps(\cdot,t)\|_{(W^{2,2}_N(\Omega))^\star}^2 dt \le C
	\qquad \mbox{for all } \eps\in (\eps_j)_{j\in\N}.
  \ee
\end{lem}
\proof
  For fixed $t\in (0,T)$ and $\psi\in W^{2,2}_N(\Omega)$, using (\ref{0eps}), (\ref{S}) and (\ref{f})
  together with the H\"older inequality, we see that since
  $\frac{\partial\psi}{\partial\nu}=0$ on $\pO$,
  \bas
	\bigg| \io \partial_t \neps (\cdot,t)\psi \bigg|
	&=& \bigg| \io \neps \Delta\psi
	+ \io \neps (S(x,\neps,\ceps)\cdot\nabla\ceps) \cdot\nabla\psi
	+ \io \neps\ueps \cdot\nabla\psi
	+ \io f(x,\neps,\ceps)\psi \bigg| \\
	&\le& \|\neps\|_{L^2(\Omega)} \|\Delta\psi\|_{L^2(\Omega)}
	+ K_S \|\neps\|_{L^\infty(\Omega)} \|\nabla\ceps\|_{L^q(\Omega)} \|\nabla\psi\|_{L^\frac{q}{q-1}(\Omega)} \\
	& &  + \|\neps\|_{L^\infty(\Omega)} \|\ueps\|_{L^r(\Omega)} \|\nabla\psi\|_{L^\frac{r}{r-1}(\Omega)}
	+ K_f \|\neps+1\|_{L^2(\Omega)} \|\psi\|_{L^2(\Omega)}
  \eas
  As the inequalities $q>N$ and $r>N$ warrant that $W^{2,2}_N(\Omega)$ is continuously embedded into both
  $W^{1,\frac{q}{q-1}}(\Omega)$ and $W^{1,\frac{r}{r-1}}(\Omega)$, this implies the existence of $C_1>0$ such that
  for all $t\in (0,T)$ and $\eps\in (\eps_j)_{j\in\N}$,
  \bas
	\|\partial_t \neps(\cdot,t)\|_{(W^{2,2}_N(\Omega))^\star}
	\le C_1 \cdot \Big\{ \|\neps\|_{L^\infty(\Omega)} +1 \Big\} \cdot
	\Big\{ \|\nabla\ceps\|_{L^q(\Omega)} + \|\ueps\|_{L^r(\Omega)} + 1 \Big\}
  \eas
  so that (\ref{36.1}) becomes a consequence of Lemma \ref{lem32} when combined with (\ref{47.1}) and (\ref{47.2}) due to
  the fact that the exponent therein satisfies $\lambda>2$.
\qed
Based on the estimates collected so far,
we can now extract a subsequence and identify a limit triple $(n,c,u)$ as follows.
\begin{lem}\label{lem38}
  Suppose that the assumptions of Theorem \ref{theo47} hold.
  Then there exist a subsequence $(\eps_{j_k})_{k\in\N}$ of $(\eps_j)_{j\in\N}$, a number $\theta\in (0,1)$ and functions
  \be{38.1}
	\left\{ \begin{array}{l}
	n\in C^{\theta,\frac{\theta}{2}}(\bom\times [0,T]), \\[1mm]
	c\in L^2((0,T);W^{1,2}(\Omega)) \qquad \mbox{and} \\[1mm]
	u\in C^{\theta,\frac{\theta}{2}}(\bom\times [0,T];\R^N)
	\cap C^{2,1}(\bom\times (0,T];\R^N)
	\end{array} \right.
  \ee
  such that as $\eps=\eps_{j_k}\searrow 0$,
  \begin{eqnarray}
	& & \neps\to n
	\qquad \mbox{in } C^0(\bom\times [0,T]),
	\label{38.2} \\[1mm]
	& & \neps\wto n
	\qquad \mbox{in }L^2((0,T); W^{1,2}(\Omega)),
	\label{38.22} \\[1mm]
	& & \ceps\wto c
	\qquad \mbox{in } L^2((0,T);W^{1,2}(\Omega))
	\qquad \mbox{and}
	\label{38.3} \\[1mm]
	& & \ueps\to u
	\qquad \mbox{in } C^0(\bom\times [0,T];\R^N) \cap C^{2,1}_{loc}(\bom\times (0,T];\R^N),
	\label{38.4}
  \eea
  and that moreover
  \be{38.5}
	\partial_t \neps \wto n_t
	\qquad \mbox{in } L^2((0,T);(W^{2,2}_N(\Omega))^\star).
  \ee
\end{lem}
\proof
  By means of a standard subsequence extraction procedure inter alia relying on the Arzel\'a-Ascoli theorem, this can
  readily be derived from Lemma \ref{lem34}, Lemma \ref{lem10},  Lemma \ref{lem48},  Lemma \ref{lem33}, Lemma \ref{lem37}
  and Lemma \ref{lem36}.
\qed
\mysection{Solution properties of $u$}
Thanks to the favorable convergence features of both $\ueps$ itself and the quantity $\neps$ determining the
forcing term in the fluid subsystem of (\ref{0eps}), it is rather evident that the limit $u$ obtained in
Lemma \ref{lem38} indeed satisfies its respective subproblem from (\ref{0e}):
\begin{lem}\label{lem39}
  If the assumptions of Theorem \ref{theo47} hold, then the functions $n$ and $u$ gained in Lemma \ref{lem38}
  have the property that with some $P\in C^{1,0}(\Omega\times (0,T))$ we have
  \be{39.1}
	u_t + \kappa (u\cdot\nabla) u = \Delta u + \nabla P + n\nabla \phi,   \quad   \nabla \cdot u=0
	\qquad \mbox{for all $x\in\Omega$ and } t\in (0,T),
  \ee
  and that $u(x,t)=0$ for all $x\in\pO$ and $t\in (0,T)$.
\end{lem}
\proof
  In view of (\ref{0eps}) and the convergence properties in (\ref{38.2}) and (\ref{38.4}), this follows from arguments
  well-established in the theory of the Navier-Stokes equations.
\qed
\mysection{Regularity and solution properties of $c$. Strong convergence of $\ceps$}
In view of the singular limit taken when passing from (\ref{0eps}) to (\ref{0e}), it may not be surprising that
corresponding questions concerning regularity in the limit process $\ceps\to c$, as well as solution properties
of the obtained limit, are more delicate.
Indeed, for appropriately taking $\eps\searrow 0$ in nonlinear expressions involving the second solution component,
and especially in the taxis term in (\ref{0eps}), the yet weak convergence information in (\ref{38.3}) seems
insufficient.
\subsection{H\"older regularity of $c, \nabla c$ and $D^2 c$. Solution properties of $c$}\label{sect5.1}
Suitable improvement of our knowledge in this respect will form the goal of this key section, and our analysis in
this direction will be launched by the following observation on validity of the Neumann problem for
second equation in (\ref{0e}) at least in some weak sense.
\begin{lem}\label{lem40}
  Let the hypotheses from Theorem \ref{theo47} be satisfied, and let $n, c$ and $u$ be as provided by Lemma \ref{lem38}.
  Then there exists a null set $N\subset (0,T)$ such that for all $t\in (0,T)\setminus N$,
  $c(\cdot,t)\in W^{1,2}(\Omega)$ with
  \be{40.1}
	\io \nabla c\cdot\nabla\psi + \io c\psi = \io n\psi - \io (u\cdot\nabla c) \psi
	\qquad \mbox{for all } \psi\in W^{1,2}(\Omega).
  \ee
\end{lem}
\proof
  Let us first make sure that for all $\varphi\in C_0^\infty(\bom\times (0,T))$ we have
  \be{40.2}
	\int_0^T \io \nabla c\cdot\nabla\varphi + \int_0^T \io c\varphi
	= \int_0^T \io n\varphi - \int_0^T \io (u\cdot\nabla c)\varphi.
  \ee
  For the verification of this, given any such $\varphi$ we use the second equation in (\ref{0eps}) to see that
  for all $\eps\in (\eps_j)_{j\in\N}$,
  \be{40.3}
	-\eps \int_0^T \io \ceps\varphi_t
	+ \int_0^T \io \nabla\ceps\cdot\nabla\varphi
	+ \int_0^T \io \ceps \varphi
	= \int_0^T \io \neps \varphi
	- \int_0^T \io (\ueps\cdot\nabla\ceps) \varphi.
  \ee
  Since (\ref{38.3}) and (\ref{38.2}) warrant that with $(\eps_{j_k})_{k\in\N}$ as found in Lemma \ref{lem38} we have
  \bas
	\int_0^T \io \nabla\ceps\cdot\nabla\varphi \to \int_0^T \io \nabla c\cdot\nabla \varphi,
	\quad
	\int_0^T \io \ceps\varphi \to \int_0^T \io c\varphi
	\quad \mbox{and} \quad
	\eps \int_0^T \io \ceps\varphi_t \to 0
  \eas
  as well as
  \bas
	\int_0^T \io \neps\varphi \to \int_0^T \io n\varphi
  \eas
  as $\eps=\eps_{j_k} \searrow 0$,
  and since combining (\ref{38.3}) with (\ref{38.4}) yields
  \bas
	\int_0^T \io (\ueps\cdot\nabla\ceps) \varphi \to \int_0^T \io (u\cdot\nabla c) \varphi
	\qquad \mbox{as } \eps=\eps_{j_k}\searrow 0,
  \eas
  the identity (\ref{40.2}) indeed results from (\ref{40.3}).\abs
  We next rely on the separability of $W^{1,2}(\Omega)$ and a mollification argument in fixing
  $(\psi_i)_{i\in\N} \subset C^\infty(\bom)$ such that $X_0:=\{\psi_i \ | \ i\in\N\}$ is dense in $W^{1,2}(\Omega)$,
  and thereupon use that all the functionals $\zeta_i^{(\iota)}$, $i\in\N, \iota \in \{1,2,3\}$, defined by
  \bas
	\zeta_i^{(1)}(t):=\io \nabla c(\cdot,t)\cdot\nabla \psi_i,
	\quad
	\zeta_i^{(2)}(t):=\io c(\cdot,t)\psi_i
	\quad \mbox{and} \quad
	\zeta_i^{(3)}(t):=\io (u(\cdot,t)\cdot\nabla c(\cdot,t))\psi_i
  \eas
  for $t\in (0,T)$ and $i\in\N$, belong to $L^1((0,T))$.
  Therefore, namely, for each $i\in\N$ we can fix a null set $\ns_i\subset (0,T)$ such that any $t\in (0,T)\setminus \ns_i$
  is a Lebesgue point of $\zeta_i^{(\iota)}$ for $\iota\in \{1,2,3\}$, whence letting
  $\ns:=\Big(\bigcup_{i\in\N} \ns_i\Big) \cup \{t\in (0,T) \ | \ c(\cdot,t)\not\in W^{1,2}(\Omega) \}$ we have found a null
  set $\ns \subset (0,T)$ such that $(0,T)\setminus \ns$ exclusively contains common Lebesgue points of all
  $\zeta_i^{(\iota)}$ for $i\in\N$ and $\iota\in \{1,2,3\}$, and such that moreover $c(\cdot,t) \in W^{1,2}(\Omega)$
  for all $t\in (0,T)\setminus \ns$.\abs
  Now for fixed $t_0\in (0,T)\setminus \ns$ and $h\in (0,T-t_0)$ we choose $(\chi_l)_{l\in\N} \subset C_0^\infty((0,T))$
  such that
  \be{40.4}
	\chi_l \wsto \chi_{(t_0,t_0+h)}
	\quad \mbox{in } L^\infty((0,T))
	\qquad \mbox{as } l\to\infty,
  \ee
  where as usual $\chi_{(t_0,t_0+h)}$ denotes the characteristic function of the set $(t_0,t_0+h)$,
  and apply (\ref{40.2}) for fixed $l\in\N$ and $\psi\in X_0$
  to $\varphi(x,t):=\chi_l(t) \cdot \psi(x)$, $(x,t)\in \Omega\times (0,T)$, to see that
  \bas
	\int_0^T \io \chi_l \nabla c\cdot\nabla\psi
	+ \int_0^T \io \chi_l c\psi
	= \int_0^T \io \chi_l n\psi
	- \int_0^T \io \chi_l (u\cdot\nabla c) \psi
	\qquad \mbox{for all  } l\in\N,
  \eas
  by (\ref{40.4}) implying that
  \bas
	\frac{1}{h} \int_{t_0}^{t_0+h} \io \nabla c\cdot\nabla \psi
	+ \frac{1}{h} \int_{t_0}^{t_0+h} \io c\psi
	= \frac{1}{h} \int_{t_0}^{t_0+h} \io n\psi
	- \frac{1}{h} \int_{t_0}^{t_0+h} \io (u\cdot\nabla c)\psi
	\quad \mbox{for all } h\in (0,T-t_0).
  \eas
  Thanks to the Lebesgue point property of $t_0$ as well as the continuity of $n$ in $\bom\times (0,T)$ asserted by
  Lemma \ref{lem38}, we may let $h\searrow 0$ here to see that
  \bas
	\io \nabla c(\cdot,t_0)\cdot\nabla \psi
	+ \io c(\cdot,t_0)\psi
	= \io n(\cdot,t_0)\psi
	- \io (u(\cdot,t_0)\cdot\nabla c(\cdot,t_0)) \psi
	\qquad \mbox{for all } \psi\in X_0,
  \eas
  which, by density of $X_0$ in $W^{1,2}(\Omega)$, upon a further approximation argument readily entails (\ref{40.1}).
\qed
Due to our knowledge on H\"older continuity of $n$ and $u$, the identity (\ref{40.1}) can be seen to entail that
$c$ actually enjoys some further regularity properties.
\begin{lem}\label{lem41}
  Under the assumptions of Theorem \ref{theo47} and with $c$ and $\ns$ taken from Lemma \ref{lem38} and Lemma \ref{lem40},
  one can find $\theta\in (0,1)$ and $C>0$ fulfilling
  \be{41.01}
	\|c(\cdot,t)\|_{W^{1,2}(\Omega)} \le C
	\qquad \mbox{for all } t\in (0,T)\setminus \ns
  \ee
  and
  \be{41.1}
	\|c(\cdot,t)-c(\cdot,s)\|_{W^{1,2}(\Omega)} \le C |t-s|^\theta
	\qquad \mbox{for all $t\in (0,T)\setminus \ns$ and } s\in (0,T)\setminus \ns.
  \ee
  In particular, on redefining $c(\cdot,t)$ for $t\in \ns\cup \{0,T\}$ if necessary, we can achieve that
  \be{41.2}
	c\in C^\theta([0,T];W^{1,2}(\Omega)).
  \ee
\end{lem}
\proof
  We first observe that for $t\in (0,T)\setminus \ns$ we may apply (\ref{40.1}) to $\psi:=c(\cdot,t) \in W^{1,2}(\Omega)$
  to see that due to Young's inequality,
  \bas
	\io |\nabla c(\cdot,t)|^2 + \io c^2(\cdot,t)
	&=& \io n(\cdot,t) c(\cdot,t) - \io (u(\cdot,t) \cdot \nabla c(\cdot,t)) c(\cdot,t) \\
	&=& \io n(\cdot,t) c(\cdot,t) \\	
	&\le& \frac{1}{2} \io c^2(\cdot,t) + \frac{1}{2} \io n^2(\cdot,t),
  \eas
  because $\nabla\cdot u(\cdot,t) \equiv 0$ in $\Omega$. By boundedness of $n$ in $\Omega\times (0,T)$, as implied
  by Lemma \ref{lem38}, this directly establishes (\ref{41.01}).\abs
  Next, for fixed $t\in (0,T)\setminus \ns$ and $s\in (0,T)\setminus \ns$, we let $z(x):=c(x,t)-c(x,s)$, $x\in\Omega$.
  Then $z\in W^{1,2}(\Omega)$ by Lemma \ref{lem40}, whence $z$ is an admissible test function in (\ref{40.1}) evaluated
  both at $t$ and at $s$.
  Subtracting the respectively obtained identities
  \bas
	\io \nabla c(\cdot,t) \cdot\nabla z + \io c(\cdot,t) z = \io n(\cdot,t) z - \io (u(\cdot,t)\cdot\nabla c(\cdot,t)) z
  \eas
  and
  \bas
	\io \nabla c(\cdot,s) \cdot\nabla z + \io c(\cdot,s) z = \io n(\cdot,s) z - \io (u(\cdot,s)\cdot\nabla c(\cdot,s)) z,
  \eas
  we thus obtain that
  \bas
	\io |\nabla z|^2 + \io z^2
	= \io (n(\cdot,t)-n(\cdot,s)) z
	- \io \Big\{ (u(\cdot,t)-u(\cdot,s))\cdot\nabla c(\cdot,t)\Big\} \cdot z
	- \io (u(\cdot,s)\cdot \nabla z) \cdot z,
  \eas
  whence if according to Lemma \ref{lem38} and (\ref{41.01}) we let $\theta_1\in (0,1)$, $C_1>0$, $C_2>0$ and $C_3>0$ be
  such that
  \bas
	& & |n(x,\tilde t)-n(x,\tilde s)| \le C_1 |\tilde t-\tilde s|^\frac{\theta_1}{2}
	\quad \mbox{and} \quad
	|u(x,\tilde t)-u(x,\tilde s)| \le C_2 |\tilde t-\tilde s|^\frac{\theta_1}{2} \\[1mm]
	& & \hspace*{48mm}
	\qquad \mbox{for all $x\in\Omega$, $\tilde t\in (0,T)$ and $\tilde s\in (0,T)$}
  \eas
  as well as
  \bas
	\|\nabla c(\cdot,\tilde t)\|_{L^2(\Omega)} \le C_3
	\qquad \mbox{for all } \tilde t\in (0,T) \setminus \ns,
  \eas
  then it follows that
  \bas
	\io |\nabla c(\cdot,t)-\nabla c(\cdot,s)|^2 + \frac{1}{2} \io (c(\cdot,t)-c(\cdot,s))^2
	\le C_1^2 |\Omega| \cdot |t-s|^{\theta_1}
	+ C_2^2 C_3^2 |t-s|^{\theta_1}
  \eas
  and that thus (\ref{41.1}) holds. The conclusion (\ref{41.2}) thereby becomes evident.
\qed
A second stage of our bootstrap-type argument now even yields some spatial $C^{2+\theta}$ regularity information, as well
as validity of the sub-problem of (\ref{0e}) in question in the classical sense:
\begin{lem}\label{lem42}
  Suppose that the assumptions from Theorem \ref{theo47} hold, and let $n, c$ and $u$ be as in Lemma \ref{lem38}.
  Then there exist $\theta\in (0,1)$ and $C>0$ such that
  \be{42.1}
	\|c(\cdot,t)\|_{C^{2+\theta}(\bom)} \le C
	\qquad \mbox{for all } t\in (0,T).
  \ee
  Moreover,
  \be{42.2}
	-\Delta c + c = n - u\cdot\nabla c
	\qquad \mbox{for all $x\in \Omega$ and } t\in (0,T)
  \ee
  as well as
  \be{42.02}
	\frac{\partial c}{\partial\nu}=0
	\qquad \mbox{for all $x\in\pO$ and } t\in (0,T)
  \ee
  in the classical sense.
\end{lem}
\proof
  Let us first make sure that there exist $q_\star>N$ and $C_1>0$ such that with $\ns\subset (0,T)$
  taken from Lemma \ref{lem40} we have
  \be{42.3}
	\|c(\cdot,t)\|_{W^{2,q_\star}(\Omega)} \le C_1
	\qquad \mbox{for all } t\in (0,T)\setminus \ns.
  \ee
  To see this, we observe that since trivially $\frac{N}{(N-2)_+}>1$, it is possible to fix a positive integer $k_0$
  and numbers $q_1,...,q_{k_0}$ such that $q_1=2, q_{k_0}>N$ and $q_k \le q_{k+1} < \frac{N q_k}{(n-q_k)_+}$ whenever
  $1\le k <k_0$.
  Then since Lemma \ref{lem38} and Lemma \ref{lem41} warrant the existence of $\theta_1\in (0,1)$, $C_2>0$, $C_3>0$ and
  $C_4>0$ such that
  \be{42.4}
	\|n(\cdot,t)\|_{C^{\theta_1}(\bom)} \le C_2,
	\quad
	\|u(\cdot,t)\|_{C^{\theta_1}(\bom)} \le C_3
	\quad \mbox{and} \quad
	\|\nabla c(\cdot,t)\|_{L^2(\Omega)} \le C_4
	\qquad \mbox{for all } t\in (0,T),
  \ee
  it follows that $h(\cdot,t):=n(\cdot,t) - u(\cdot,t)\cdot\nabla c(\cdot,t)$, $t\in (0,T)$, has the property that
  by definition of $q_1$,
  \bas
	\|h(\cdot,t)\|_{L^{q_1}(\Omega)} \le C_5:=C_2 + C_3 C_4
	\qquad \mbox{for all } t\in (0,T).
  \eas
  As from Lemma \ref{lem40} we know that for each $t\in (0,T)\setminus \ns$, $c(\cdot,t)\in W^{1,2}(\Omega)$ is a weak
  solution, in the standard sense specified in (\ref{40.1}), of the Neumann boundary value problem for
  $-\Delta c(\cdot,t)+c(\cdot,t)=h(\cdot,t)$ in $\Omega$, elliptic estimates therefor (\cite{GT}) provide $C_6>0$ such that
  \bas
	\|c(\cdot,t)\|_{W^{2,q_1}(\Omega)} \le C_6\|h(\cdot,t)\|_{L^{q_1}(\Omega)} \le C_5 C_6
	\qquad \mbox{for all } t\in (0,T)\setminus \ns.
  \eas
  In the case $N=1$ in which $q_1>N$, this already establishes (\ref{42.3}), while if $N\ge 2$ and hence $k_0>1$,
  the inequality $q_2<\frac{N q_1}{(N-q_1)_+}$ ensures continuity of the embedding $W^{2,q_1}(\Omega) \hra W^{1,q_2}(\Omega)$,
  whence (\ref{42.4}) actually implies boundedness of $(h(\cdot,t))_{t\in (0,T)\setminus\ns}$ in $L^{q_2}(\Omega)$.\abs
  Repeating this procedure, after finitely many steps we conclude that indeed (\ref{42.3}) holds with $q_\star:=q_{k_0}$
  and some appropriately large $C_1>0$. Since $q_\star>N$, in view of the continuous embedding
  $W^{2,q_\star}(\Omega) \hra C^{1+\theta_2}(\bom)$ for any fixed $\theta_2\in (0,1-\frac{N}{q_\star})$ this entails
  boundedness of $(\nabla c(\cdot,t))_{t\in (0,T)\setminus\ns}$ in $C^{\theta_2}(\bom)$ and thus, again through
  (\ref{42.4}), of $(h(\cdot,t))_{t\in (0,T)\setminus\ns}$ in $C^{\theta_3}(\bom)$ for some $\theta_3\in (0,1)$.
  Now elliptic Schauder theory (\cite{GT}) applies to provide $C_7>0$ such that
  \bas
	\|c(\cdot,t)\|_{C^{2+\theta_3}(\bom)} \le C_7
	\qquad \mbox{for all } t\in (0,T)\setminus\ns,
  \eas
  which due to the time continuity property expressed in (\ref{41.2}) clearly extends so as to remain valid actually for
  all $t\in (0,T)$. This clearly entails (\ref{42.1}) and, as a consequence of (\ref{40.1}), also (\ref{42.2})  and thus (\ref{42.02}).
\qed
By straightforward interpolation, combining the latter two lemmata finally provides H\"older continuity also in time
of $c, \nabla c$ and $D^2 c$.
\begin{lem}\label{lem43}
  Under the assumptions from Theorem \ref{theo47}, there exist $\theta\in (0,1)$ and $C>0$ such that
  \be{43.1}
	\|c(\cdot,t)-c(\cdot,s)\|_{C^{2+\theta}(\bom)} \le C|t-s|^\theta
	\qquad \mbox{for all $t\in (0,T)$ and } s\in (0,T),
  \ee
  where $c$ is taken from Lemma \ref{lem38}.
\end{lem}
\proof
  In line with Lemma \ref{lem41}, let us fix $\theta_1\in (0,1)$ such that $c\in C^{\theta_1}([0,T];W^{1,2}(\Omega))$,
  and thereafter choose any $\theta_2\in (0,\theta_1)$.
  Then by straightforward interpolation, we can find $a\in (0,1)$ and $C_1>0$ such that
  \bas
	\|\psi\|_{C^{2+\theta_2}(\bom)} \le C_1 \|\psi\|_{C^{2+\theta_1}(\bom)}^a \|\psi\|_{W^{1,2}(\Omega)}^{1-a}
	\qquad \mbox{for all } \psi \in C^{2+\theta_1}(\bom).
  \eas
  Therefore,
  \bas
	\|c(\cdot,t)-c(\cdot,s)\|_{C^{2+\theta_2}(\bom)}
	\le C_1 \cdot \Big\{ \|c(\cdot,t)\|_{C^{2+\theta_1}(\bom)} + \|c(\cdot,s)\|_{C^{2+\theta_1}(\bom)} \Big\}^a
	\cdot \|c(\cdot,t)-c(\cdot,s)\|_{W^{1,2}(\Omega)}^{1-a}
  \eas
  for all $t\in (0,T)$ and $s\in (0,T)$,
  so that the claim readily results from Lemma \ref{lem42} and Lemma \ref{lem41} if we let $\theta:=(1-a)\theta_1$
  and take $C>0$ appropriately large.
\qed
\subsection{Regularity of $c_t$. Strong convergence properties of $\ceps$ and $\nabla\ceps$}\label{sect5.2}
Now an observation crucial for our derivation of strong convergence properties of $\ceps$ is contained
in the following statement on local $L^2$ integrability of $c_t$ away from $t=0$.
Using the validity of the elliptic subproblem in (\ref{0e}) for $c$ as a starting point, besides relying on
the local boundedness of $u_t$ in $\bom\times (0,T]$ our argument essentially utilizes the time regularity
information on $n_t$ provided by Lemma \ref{lem38} through Lemma \ref{lem36}.
\begin{lem}\label{lem44}
  Under the assumptions from Theorem \ref{theo47}, the function $c$ obtained in Lemma \ref{lem38} satisfies
  \be{44.1}
	c_t \in L^2_{loc}(\bom\times (0,T]).
  \ee
\end{lem}
\proof
  We pick $\tau\in (0,T)$ and $h_0\in (0,T-\tau)$, and for $h\in (0,h_0)$ we let
  \bas
	z_h(x,t):=\frac{c(x,t+h)-c(x,t)}{h},
	\qquad x\in\Omega, \ t\in (\tau,T-h_0).
  \eas
  Then using Lemma \ref{lem42} we see that for each fixed $t\in (\tau,T-h_0)$, $z_h(\cdot,t) \in C^2(\bom)$ is a classical
  solution of the Neumann boundary value problem for
  \bea{44.2}
	-\Delta z_h(\cdot,t) + z_h(\cdot,t)
	&=& \frac{1}{h} \cdot \Big\{ n(\cdot,t+h) - u(\cdot,t+h)\cdot\nabla c(\cdot,t+h)\Big\}
	- \frac{1}{h} \cdot \Big\{ n(\cdot,t) - u(\cdot,t)\cdot\nabla c(\cdot,t)\Big\} \nn\\
	&=& g_h(\cdot,t) - u(\cdot,t)\cdot\nabla z_h(\cdot,t)
	\qquad \mbox{in } \Omega,
  \eea
  where for $h\in (0,h_0)$,
  \bas
	g_h(x,t):=\frac{n(x,t+h)-n(x,t)}{h} - \frac{u(x,t+h)-u(x,t)}{h} \cdot \nabla c(x,t+h),
	\qquad x\in\Omega, \ t\in (\tau,T-h_0).
  \eas
  Next, for $h\in (0,h_0)$ and $t\in (0,T-h_0)$ we furthermore let $\zho(\cdot,t)$ and $\zht(\cdot,t)$ denote the classical
  solutions of
  \be{44.3}
	\left\{ \begin{array}{ll}
	-\Delta \zho(\cdot,t) + \zho(\cdot,t) = g_h(\cdot,t),
	\qquad & x\in\Omega, \\[1mm]
	\frac{\partial \zho}{\partial\nu}=0,
	& x\in\pO,
	\end{array} \right.
  \ee
  and
  \be{44.4}
	\left\{ \begin{array}{ll}
	-\Delta \zht(\cdot,t) + \zht(\cdot,t) = - u(\cdot,t)\cdot \nabla z_h(\cdot,t),
	\qquad & x\in\Omega, \\[1mm]
	\frac{\partial \zht}{\partial\nu}=0,
	& x\in\pO,
	\end{array} \right.
  \ee
  noting that their existence in the space $C^2(\bom)$ is asserted by standard elliptic theory (\cite{GT}) due to the
  fact that both $g_h(\cdot,t)$ and $u(\cdot,t)\cdot\nabla z_h(\cdot,t)$ are H\"older continuous in $\bom$ by
  Lemma \ref{lem38}, Lemma \ref{lem42} and the inclusion $z_h(\cdot,t) \in C^2(\bom)$.\abs
  Then according to a uniqueness property of classical solutions to the Neumann problem associated with the inhomogeneous
  Helmholtz equation in (\ref{44.2}), it follows that
  \be{44.44}
	z_h(\cdot,t)=\zho(\cdot,t)+\zht(\cdot,t)
	\qquad \mbox{for all } t\in (\tau,T-t_0),
  \ee
  and in order to successively derive $L^2$ bounds for $\zho$ and then for $\zht$, we first rewrite (\ref{44.3}) in the form
  $\zho(\cdot,t):=B^{-1} g_h(\cdot,t)$, where $B$ denotes the realization of $-\Delta +1$
  in $W^{2,2}_N(\Omega)$.
  Then since $B^{-1}$ obviously is nonexpansive on $L^2(\Omega)$, we can estimate
  \bas
	\|\zho(\cdot,t)\|_{L^2(\Omega)}
	&\le& \bigg\| B^{-1} \frac{n(\cdot,t+h)-n(\cdot,t)}{h}\bigg\|_{L^2(\Omega)}
	+ \bigg\| B^{-1} \Big\{ \frac{u(\cdot,t+h)-u(\cdot,t)}{h} \cdot\nabla c(\cdot,t+h)\Big\} \bigg\|_{L^2(\Omega)} \\
	&\le& \bigg\| B^{-1} \frac{n(\cdot,t+h)-n(\cdot,t)}{h}\bigg\|_{L^2(\Omega)}
	+ \bigg\| \frac{u(\cdot,t+h)-u(\cdot,t)}{h} \cdot\nabla c(\cdot,t+h) \bigg\|_{L^2(\Omega)} \\
	&\le& \bigg\| B^{-1} \frac{n(\cdot,t+h)-n(\cdot,t)}{h}\bigg\|_{L^2(\Omega)}
	+ \bigg\| \frac{u(\cdot,t+h)-u(\cdot,t)}{h} \bigg\|_{L^\infty(\Omega)} \|\nabla c(\cdot,t+h)\|_{L^2(\Omega)} \\
	&=& \bigg\| \frac{1}{h} \int_t^{t+h} B^{-1} n_t(\cdot,s) ds \bigg\|_{L^2(\Omega)}
	+ \bigg\| \frac{1}{h} \int_t^{t+h} u_t(\cdot,s) ds \bigg\|_{L^\infty(\Omega)}
		\cdot \|\nabla c(\cdot,t+h)\|_{L^2(\Omega)} \\[2mm]
	& & \hspace*{55mm} \mbox{for all $t\in (\tau,T-h_0)$ and any } h\in (0,h_0).
  \eas
  Since $u_t$ is bounded in $\Omega\times (\tau,T)$ by Lemma \ref{lem38}, and since Lemma \ref{lem41} implies boundedness
  of $(0,T) \ni t \mapsto \|\nabla c(\cdot,t)\|_{L^2(\Omega)}$, we thus obtain $C_1=C_1(\tau)>0$ such that for all
  $h\in (0,h_0)$,
  \bas
	\|\zho(\cdot,t)\|_{L^2(\Omega)}^2
	\le 2\bigg\| \frac{1}{h} \int_t^{t+h} B^{-1} n_t(\cdot,s) ds \bigg\|_{L^2(\Omega)}^2 + C_1
	\qquad \mbox{for all } t\in (\tau,T-h_0).
  \eas
  Therefore, by integration using the Fubini theorem as well as the Cauchy-Schwarz inequality,
  \bas
	\int_\tau^{T-h_0} \|\zho(\cdot,t)\|_{L^2(\Omega)}^2 dt
   	&\le& C_1 T + \frac{2}{h^2} \int_\tau^{T-h_0} \bigg\| \int_t^{t+h} B^{-1} n_t(\cdot,s) ds \bigg\|_{L^2(\Omega)}^2 dt
		\\
	&\le& C_1 T + \frac{2}{h} \int_\tau^{T-h_0} \int_t^{t+h} \|B^{-1} n_t(\cdot,s)\|_{L^2(\Omega)}^2 dsdt \\
	&=& C_1 T + \frac{2}{h} \int_\tau^{\tau+h} \int_\tau^s \|B^{-1} n_t(\cdot,s)\|_{L^2(\Omega)}^2 dtds \\
	& & + \frac{2}{h} \int_{\tau+h}^{T-h_0} \int_{s-h}^s \|B^{-1} n_t(\cdot,s)\|_{L^2(\Omega)}^2 dtds \\
	& & + \frac{2}{h} \int_{T-h_0}^{T-h_0+h} \int_{s-h}^{T-h_0} \|B^{-1} n_t(\cdot,s)\|_{L^2(\Omega)}^2 dtds \\
	&\le& C_1 T + 2 \int_\tau^T \| B^{-1} n_t(\cdot,s)\|_{L^2(\Omega)}^2 ds
	\qquad \mbox{for all } h\in (0, T-h_0-\tau).
  \eas
  Now since standard elliptic regularity theory (\cite{wloka}) ensures that $B^{-1}$ maps $L^2(\Omega)$ continuously
  into $W^{2,2}_N(\Omega)$, and that thus with some $C_2>0$ we have $\|B^{-1}\psi\|_{L^2(\Omega)} \le
  C_2 \|\psi\|_{(W^{2,2}_N(\Omega))^\star}$ for all $\psi\in (W^{2,2}_N(\Omega))^\star$, this entails that
  \be{44.5}
	\int_\tau^{T-h_0} \io (\zho)^2 	
	\le C_3=C_3(\tau):=C_1 T + 2C_2^2 \int_\tau^T \|n_t(\cdot,t)\|_{(W^{2,2}_N(\Omega))^\star}^2 dt
  \ee
  for all $h\in 0, T-h_0-\tau)$  with $C_3$ being finite thanks to Lemma \ref{lem36}.\abs
  Next, in order to estimate $\zht$, we test (\ref{44.4}) against $\zht(\cdot,t)$ and recall the decomposition
  (\ref{44.44}) to infer that by solenoidality of $u$ and by Young's inequality,
  \bas
	\io |\nabla \zht(\cdot,t)|^2 + \io (\zht(\cdot,t))^2
	&=& - \io (u(\cdot,t)\cdot\nabla \zho(\cdot,t)) \zht(\cdot,t)
	- \io (u(\cdot,t)\cdot\nabla \zht(\cdot,t)) \zht(\cdot,t) \\
	&=&  \io \zho(\cdot,t) (u(\cdot,t)\cdot\nabla \zht(\cdot,t)) \\
	&\le& \io |\nabla\zht(\cdot,t)|^2
	+ \frac{1}{4} \io (\zho(\cdot,t))^2 |u(\cdot,t)|^2 \\
	&\le& \io |\nabla\zht(\cdot,t)|^2
	+ \frac{1}{4} \|u(\cdot,t)\|_{L^\infty(\Omega)}^2 \io (\zho(\cdot,t))^2 \\[2mm]
	& & \hspace*{8mm} \qquad \mbox{for all $t\in (\tau,T-h_0)$ and each }   h\in (0, T-h_0-\tau).
  \eas
  As once more relying on Lemma \ref{lem38} we can find $C_4>0$ such that $\|u(\cdot,t)\|_{L^\infty(\Omega)} \le C_4$
  for all $t\in (0,T)$, in view of (\ref{44.5}) this implies that
  \bas
	\int_\tau^{T-h_0} \io (\zht)^2
	\le \frac{C_4^2}{4} \int_\tau^{T-h_0} \io (\zho)^2
	\le \frac{C_3 C_4^2}{4}
	\qquad \mbox{for all } h\in  (0, T-h_0-\tau)
  \eas
  and that thus, again by (\ref{44.5}),
  \bas
	\int_\tau^{T-h_0} \io z_h^2 \le 2C_3 + \frac{C_3 C_4^2}{2}
	\qquad \mbox{for all } h\in  (0, T-h_0-\tau)
  \eas
  because of (\ref{44.44}).
  Consequently, there exist   $(h_i)_{i\in\N} \subset (0, T-h_0-\tau)$ and $z\in L^2(\Omega\times (\tau, T-h_0))$
  such that $h_i\to 0$ and $z_{h_i} \wto z$ in $L^2(\Omega\times (\tau,T-h_0))$ as $i\to\infty$, so that since by definition
  of distributional derivatives $z$ must coincide with $c_t$ a.e.~in $\Omega\times (\tau,T-h_0)$, (\ref{44.1})
  results from the fact that  $\tau\in (0,T)$ and $h_0\in (0, T-\tau)$  were arbitrary.
\qed
On the basis of this, we can now in fact derive some strong convergence property of $\ceps$ by
analyzing the difference $\ceps-c$ through a parabolic equation satisfied
by the latter, in which the crucial source term $\eps c_t$, namely, can appropriately be controlled using
(\ref{44.1}):
\begin{lem}\label{lem45}
  Suppose that the assumptions of Theorem \ref{theo47} are satisfied with some
  $(\eps_j)_{j\in\N} \subset (0,\infty)$, and let $(n,c,u,P)$ and $(\eps_{j_k})_{k\in\N}$ be as provided by Lemma \ref{lem38}.
  Then
  \be{45.1}
	\ceps\to c
	\qquad \mbox{in } L^\infty_{loc}((0,T];L^2(\Omega))
  \ee
  and
  \be{45.2}
	\nabla\ceps\to \nabla c
	\qquad \mbox{in } L^2_{loc}(\bom\times (0,T])
  \ee
  as $\eps=\eps_{j_k}\searrow 0$.
\end{lem}
\proof
  For $\eps\in (\eps_j)_{j\in\N}$, we let
  \bas
	\zeps(x,t):=\ceps(x,t) - c(x,t),
	\qquad (x,t)\in\Omega\times (0,T),
  \eas
  and
  \bas
	\yeps(t):=\io \zeps^2(x,t) dx,
	\qquad t\in (0,T).
  \eas
  Then since $\zeps\in L^\infty(\Omega\times (0,T))$ and
  $\partial_t \zeps=\partial_t \ceps-c_t\in L^2_{loc}((0,T];L^2(\Omega))$ by Lemma \ref{lem01} with $q=2$ and Lemma \ref{lem44},
  it follows from a standard argument that $\yeps$ belongs to $W^{1,2}_{loc}((0,T])$ and is therefore locally absolutely
  continuous in $(0,T]$ with
  \bas
	\yeps'(t) = 2\io \zeps(\cdot,t) \partial_t \zeps(\cdot,t)
	\qquad \mbox{for a.e.~} t\in (0,T).
  \eas
  As herein by (\ref{0eps}) and Lemma \ref{lem42},
  \bas
	\eps\partial_t \zeps
	&=& \Delta \ceps - \ceps + \neps - \ueps\cdot\nabla\ceps - \eps c_t \\
	&=& \Delta \zeps + \Delta c - \zeps - c  + \neps - \ueps\cdot\nabla \ceps - \eps c_t \\
	&=& \Delta \zeps - \zeps + (\neps-n) - (\ueps-u) \cdot \nabla c - \ueps\cdot\nabla\zeps - \eps c_t
	\qquad \mbox{a.e.~in $\Omega$ for a.e.~} t\in (0,T),
  \eas
  on integrating by parts and using Young's inequality we obtain that
  \bas
	& & \hspace*{-20mm}
	\frac{\eps}{2} \yeps'(t)
	+ \io |\nabla\zeps|^2 + \io \zeps^2 \nn\\
	&=& \io (\neps-n) \zeps
	- \io \Big\{ (\ueps-u)\cdot\nabla c\Big\} \zeps
	- \io (\ueps\cdot\nabla\zeps) \zeps
	- \eps \io c_t \zeps \nn\\
	&\le& \frac{3}{4} \io \zeps^2
	+ \io (\neps-n)^2
	+ \io |\ueps-u|^2 |\nabla c|^2
	+ \eps^2 \io c_t^2
	\qquad \mbox{for a.e.~} t\in (0,T),
  \eas
  because $\nabla\cdot \ueps \equiv 0$ a.e.~in $\Omega\times (0,T)$. Hence
  \be{45.3}
	\frac{\eps}{2} \yeps'(t)
	+ \frac{1}{4} \yeps(t) + \io |\nabla\zeps|^2
	\le \heps + \eps^2 \io c_t^2(\cdot,t)
	\qquad \mbox{for a.e.~} t\in (0,T),
  \ee
  where according to Lemma \ref{lem38},
  \bas
	\heps:=|\Omega| \cdot \|\neps-n\|_{L^\infty(\Omega\times (0,T))}^2
	+ \|\ueps-u\|_{L^\infty(\Omega\times (0,T))}^2 \cdot \|\nabla c\|_{L^\infty((0,T);L^2(\Omega))}^2
  \eas
  satisfies
  \be{45.4}
	\heps \to 0
	\qquad \mbox{as } \eps=\eps_{j_k}\searrow 0.
  \ee
  To make sure that in conjunction with Lemma \ref{lem44} and Lemma \ref{lem48} this implies that for each $\tau\in (0,T)$
  we have
  \be{45.5}
	\zeps \to 0
	\quad \mbox{in } L^\infty((\tau,T);L^2(\Omega))
	\qquad \mbox{as } \eps=\eps_{j_k}\searrow 0,
  \ee
  given any such $\tau$ and an arbitrary $\eta>0$ we use (\ref{45.4}) and Lemma \ref{lem44} along with the boundedness
  of $(\yeps)_{\eps\in (\eps_{j_k})_{k\in\N}}$ in $L^\infty((0,T))$, as entailed by Lemma \ref{lem48} and Lemma \ref{lem42},
  to fix $\eps_0>0$ small enough such that whenever $\eps\in (\eps_{j_k})_{k\in\N}$ is such that $\eps<\eps_0$, we have
  \be{45.6}
	4\heps \le \frac{\eta}{3}
  \ee
  and
  \be{45.7}
	2\eps \int_{\frac{\tau}{2}}^T \io c_t^2 \le \frac{\eta}{3}
  \ee
  as well as
  \be{45.8}
	\yeps\Big(\frac{\tau}{2}\Big) \cdot e^{-\frac{\tau}{4\eps}} \le \frac{\eta}{3}.
  \ee
  Then from (\ref{45.3}) we infer on dropping the nonnegative last summand on the left that the absolutely continuous function
  $[\frac{\tau}{2}, T] \ni t \mapsto e^{\frac{1}{2\eps} (t-\frac{\tau}{2})} \yeps(t)$ satisfies
  \bas
	\frac{d}{dt} \bigg\{ e^{\frac{1}{2\eps}(t-\frac{\tau}{2})} \yeps(t) \bigg\}
	\le e^{\frac{1}{2\eps}(t-\frac{\tau}{2})}  \bigg\{   \frac{2\heps}{\eps} + 2\eps \io c_t^2(\cdot,t) \bigg\}
	\qquad \mbox{for a.e.~} t\in \Big(\frac{\tau}{2},T\Big),
  \eas
  and therefore we obtain using (\ref{45.8}), (\ref{45.6}) and (\ref{45.7}) that
  \bas
	\yeps(t)
	&\le& \yeps\Big(\frac{\tau}{2}\Big) \cdot e^{-\frac{1}{2\eps}(t-\frac{\tau}{2})}
	+ \frac{2\heps}{\eps} \int_{\frac{\tau}{2}}^t e^{-\frac{1}{2\eps}(t-s)} ds
	+ 2\eps \int_{\frac{\tau}{2}}^t e^{-\frac{1}{2\eps}(t-s)} \cdot \io c_t^2(\cdot,s) ds \\
	&\le& \yeps\Big(\frac{\tau}{2}\Big) \cdot e^{-\frac{1}{2\eps}(t-\frac{\tau}{2})}
	+ 4\heps \cdot \Big\{ 1- e^{-\frac{1}{2\eps}(t-\frac{\tau}{2})} \Big\}
	+ 2\eps \int_{\frac{\tau}{2}}^t \io c_t^2 \\[1mm]
	&\le& \frac{\eta}{3} + \frac{\eta}{3} + \frac{\eta}{3}=\eta
	\qquad \mbox{for all }  t\in \big(\tau,T\big).
  \eas
  Having thereby verified (\ref{45.5}), going back to (\ref{45.3})we see upon direct integration therein that for each
  $\tau\in (0,T)$ we moreover have
  \bas
	\int_\tau^T \io |\nabla\zeps|^2
	\le \frac{\eps}{2} \yeps(\tau) + \heps\cdot (T-\tau) + \eps^2 \int_\tau^T \io c_t^2
	\qquad \mbox{for all } \eps\in (\eps_{j_k})_{k\in\N},
  \eas
  so that again by means of (\ref{45.4}), Lemma \ref{lem44} and the boundedness of $(\yeps)_{\eps\in (\eps_{j_k})_{k\in\N}}$
  in $L^\infty((0,T))$ we infer that for any such $\tau$,
  \bas
	\zeps \to 0
	\quad \mbox{in } L^2((\tau,T);W^{1,2}(\Omega))
	\qquad \mbox{as } \eps=\eps_{j_k}\searrow 0.
  \eas
  Together with (\ref{45.5}), this shows that both (\ref{45.1}) and (\ref{45.2}) hold.
\qed
\subsection{A bound for $\nabla \ceps$ in $L^\infty((0,T);L^{\hatq}(\Omega))$ for arbitrarily large $\hatq$}
Let us conclude this section by providing some additional integrability information on the signal gradient on the basis
of the differential inequality from Lemma \ref{lem01}.\abs
Our reasoning will involve the following elementary interpolation inequality.
\begin{lem}\label{lem49}
  Let $q\ge 2$. Then for all $\varphi\in C^2(\bom)$ such that $\varphi \cdot \frac{\partial\varphi}{\partial\nu}=0$ on $\pO$,
  we have
  \be{49.1}
	\io |\nabla\varphi|^{q}
	\le (\sqrt{N}+q-2)^{\frac{2q}{q+2}} \cdot \bigg\{ \io |\nabla\varphi|^{q-2} |D^2 \varphi|^2 \bigg\}^\frac{q}{q+2}
	\cdot \bigg\{ \io |\varphi|^q \bigg\}^\frac{2}{q+2}.
  \ee
\end{lem}
\proof
  We integrate by parts and use the H\"older inequality to see that
  \bas
	\io |\nabla\varphi|^q
	&=& \io |\nabla\varphi|^{q-2} \nabla\varphi\cdot\nabla\varphi \\
	&=& - \io \varphi |\nabla \varphi|^{q-2} \Delta\varphi
	- (q-2) \io \varphi |\nabla \varphi|^{q-4} \nabla\varphi\cdot (D^2\varphi \cdot\nabla\varphi) \\
	&\le& (\sqrt{N}+q-2) \io |\varphi| \cdot |\nabla\varphi|^{q-2} \cdot |D^2\varphi| \\
	&\le& (\sqrt{N}+q-2) \bigg\{ \io |\varphi|^q \bigg\}^\frac{1}{q} \cdot
	\bigg\{ \io |\nabla \varphi|^q \bigg\}^\frac{q-2}{2q} \cdot
	\bigg\{ \io |\nabla\varphi|^{q-2} |D^2\varphi|^2 \bigg\}^\frac{1}{2},
  \eas
  from which (\ref{49.1}) can readily be derived.
\qed	
Indeed, we can thereby achieve the following.
\begin{lem}\label{lem50}
  Assume the hypotheses from Theorem \ref{theo47}, and
  let $\hatq>  \max \{N,2\}$. Then there exists $C=C(\hatq)>0$ such that for all $\eps\in (\eps_j)_{j\in\N}$,
  \be{50.1}
	\io |\nabla\ceps(\cdot,t)|^{\hatq} \le C
	\qquad \mbox{for all } \,  t\in (0,T).
  \ee
  In particular,
  \be{50.2}
	\sup_{\eps \in (\eps_j)_{j\in\N}} \|\ceps\|_{L^\infty(\Omega\times (0,T))} <\infty.
  \ee
\end{lem}
\proof
  Since $(\neps)_{\eps\in (\eps_j)_{j\in\N}}$ and $(\ueps)_{\eps\in (\eps_j)_{j\in\N}}$
  are bounded in $L^\infty(\Omega\times (0,T))$
  and in $L^\infty(\Omega\times (0,T);\R^N)$ according to Lemma \ref{lem32} and Lemma \ref{lem33}, respectively,
  from Lemma \ref{lem01} we infer the existence of $C_1=C_1(\hatq)>0$ such that for all $\eps\in (\eps_j)_{j\in\N}$,
  \bea{50.3}
	\frac{\eps}{\hatq} \frac{d}{dt} \io |\nabla\ceps|^{\hatq}
	+ \frac{1}{4} \io |\nabla \ceps|^{\hatq-2} |D^2\ceps|^2 + \io |\nabla\ceps|^{\hatq}
	\le C_1 + C_1 \io |\nabla \ceps|^{\hatq}
	\quad  \, \, \mbox{for all } \,  t\in (0,T).
  \eea
  Here we may use that Lemma \ref{lem48} warrants boundedness of $(\ceps)_{\eps\in (\eps_j)_{j\in\N}}$ in
  $L^\infty((0,T);L^{\hatq}(\Omega))$ to see that as a consequence of Lemma \ref{lem49} and Young's inequality, with
  some positive constants $C_i=C_i(\hatq)$, $i\in \{2,3,4\}$, we have
  \bas
	C_1 \io |\nabla\ceps|^{\hatq}
	&\le& C_2 \cdot \bigg\{ \io |\nabla\ceps|^{\hatq-2} |D^2\ceps|^2 \bigg\}^\frac{\hatq}{\hatq+2} \cdot
	\bigg\{ \io \ceps^{\hatq} \bigg\}^\frac{2}{q+2} \\
	&\le& C_3 \cdot \bigg\{ \io |\nabla\ceps|^{\hatq-2} |D^2\ceps|^2 \bigg\}^\frac{\hatq}{\hatq+2} \\
	&\le& \frac{1}{4} \io |\nabla\ceps|^{\hatq-2} |D^2\ceps|^2 + C_4
	\qquad \mbox{for all $t\in (0,T)$ and any } \eps\in (\eps_j)_{j\in\N},
  \eas
  so that (\ref{50.3}) entails the inequality
  \bas
	\frac{\eps}{\hatq} \frac{d}{dt} \io |\nabla\ceps|^{\hatq}
	+ \io |\nabla\ceps|^{\hatq}
	\le C_1 + C_4
	\qquad \mbox{for all $t\in (0,T)$ and } \eps\in (\eps_j)_{j\in\N},
  \eas
  from which (\ref{50.1}) directly follows.
  As $\hatq>N$ and hence $W^{1,\hatq}(\Omega) \hra L^\infty(\Omega)$, again due to the boundedness property from
  Lemma \ref{lem48} this in turn implies (\ref{50.2}).
\qed
\mysection{Solution properties of $n$. Proof of Theorem \ref{theo47}}\label{sect6}
Now knowing that $\ceps\to c$ also in the pointwise sense, we can readily pass to the limit also in the first
equation in (\ref{0eps}).
\begin{lem}\label{lem46}
  Suppose that the assumptions of Theorem \ref{theo47} are satisfied with some
  $(\eps_j)_{j\in\N} \subset (0,\infty)$, and let $(n,c,u,P)$ be as provided by Lemma \ref{lem38}.
  Then in the classical pointwise sense we have
  \be{46.1}
	n_t + u\cdot\nabla n = \Delta n - \nabla\cdot (nS(x,n,c)\cdot\nabla c) + f(x,n,c)
	\qquad \mbox{for all $x\in\Omega$ and } t\in (0,T)
  \ee
  as well as
  \be{46.01}
	(\nabla n - nS(x,n,c)\cdot\nabla c)\cdot \nu=0
	\qquad \mbox{for all $x\in\pO$ and } t\in (0,T).
  \ee
\end{lem}
\proof
  Let us first make sure that for arbitrary $\varphi\in C_0^\infty(\bom\times [0,T))$,
  \bea{46.2}
	- \int_0^T \io n \varphi_t - \io n_0 \varphi(\cdot,0)
	&=& - \int_0^T \io \nabla n \cdot\nabla\varphi
	+ \int_0^T \io n (S(x,n,c)\cdot\nabla c)\cdot\nabla\varphi \nn\\
	& & + \int_0^T \io nu\cdot\nabla\varphi
	+ \int_0^T \io f(x,n,c)\varphi.
  \eea
  To see this, given any such $\varphi$ and $\eps\in (\eps_j){j\in\N}$ we use (\ref{0eps}) to find that
  \bea{46.3}
	- \int_0^T \io \neps \varphi_t - \io n_0 \varphi(\cdot,0)
	&=& - \int_0^T \io \nabla \neps \cdot\nabla\varphi
	+ \int_0^T \io \neps (S(x,\neps,\ceps)\cdot\nabla \ceps)\cdot\nabla\varphi \nn\\
	& & + \int_0^T \io \neps\ueps\cdot\nabla\varphi
	+ \int_0^T \io f(x,\neps,\ceps)\varphi,
  \eea
  where by (\ref{38.2}), (\ref{38.22}) and (\ref{38.4}), clearly
  \bas
	& & \int_0^T \io \neps\varphi_t \to \int_0^T \io n\varphi_t,
	\qquad \int_0^T \io \nabla\neps\cdot\nabla\varphi \to \int_0^T \io \nabla n\cdot\nabla \varphi
	\qquad \mbox{and} \\
	& & \int_0^T \io \neps\ueps \cdot\nabla\varphi \to \int_0^T \io nu\cdot\nabla\varphi
  \eas
  as $\eps=\eps_{j_k} \searrow 0$, where $(\eps_{j_k})_{k\in\N}$ is as provided by Lemma \ref{lem38}.
  Apart from this, thanks to Lemma \ref{lem45} we know that $\ceps \to c$ a.e.~in $\Omega\times (0,T)$ and hence,
  again by (\ref{38.2}), that also $S(\cdot,\neps,\ceps) \to S(\cdot,n,c)$ and $f(\cdot,\neps,\ceps) \to f(\cdot,n,c)$
  a.e.~in $\Omega\times (0,T)$.
  Since combining (\ref{S}) and (\ref{f}) with (\ref{38.2}) and (\ref{50.2}) moreover ensures boundedness of
  $(S(\cdot,\neps,\ceps))_{\eps\in (\eps_{j_k})_{k\in\N}}$ in $L^\infty(\Omega\times (0,T);\R^{N\times N})$
  and of $(f(\cdot,\neps,\ceps))_{\eps\in (\eps_{j_k})_{k\in\N}}$ in $L^\infty(\Omega\times (0,T))$, by means of the
  dominated convergence theorem and a well-known argument (\cite[Lemma A.4]{win_ct_rot_SIMA}) we conclude that as
  $\eps=\eps_{j_k}\searrow 0$, not only $f(\cdot,\neps, \ceps) \to f(\cdot,n,c)$ in $L^1(\Omega\times (0,T))$ but also
  \bas
	\neps S(\cdot,\neps,\ceps) \to nS(\cdot,n,c)
	\qquad \mbox{in } L^2(\Omega\times (0,T);\R^{N\times N}).
  \eas
  In conjunction with (\ref{38.3}), these properties warrant that
  \bas
	\int_0^T \io f(x,\neps,\ceps) \varphi \to \int_0^T \io f(x,n,c)\varphi
  \eas
  and
  \bas
	\int_0^T \io \neps (S(x,\neps,\ceps)\cdot\nabla\ceps) \cdot \nabla\varphi
	\to \int_0^T \io n (S(x,n,c)\cdot\nabla c)\cdot\nabla \varphi
  \eas
  as $\eps=\eps_{j_k}\searrow 0$, so that (\ref{46.2}) becomes a consequence of (\ref{46.3}).\abs
  Thus knowing that $n\in L^2((0,T);W^{1,2}(\Omega))$ forms a weak solution of the initial-boundary value problem
  (\ref{46.1})-(\ref{46.01}) in the standard generalized sense from e.g.~\cite{LSU}, we may invoke classical results from
  parabolic regularity theory to conclude from the H\"older continuity of $n, c,\nabla c, D^2 c$ and $u$ in
  $\bom\times [0,T]$, as stated by Lemma \ref{lem38} and Lemma \ref{lem43}, that firstly
  $n\in C^{1+\theta_1, \frac{1+\theta_1}{2}}(\bom\times (0,T])$ for some $\theta_1\in (0,1)$ (\cite{lieberman}),
  and that secondly, as a consequence thereof, due to this additional information on H\"older continuity of
  $\nabla n$ we can find $\theta_2\in (0,1)$ such that we even have
  $n\in C^{2+\theta_2,1+\frac{\theta_2}{2}}(\bom\times (0,T])$,
  and that therefore (\ref{46.2}) warrants validity of (\ref{46.1}) and (\ref{46.01}) in the classical pointwise sense
  through a standard variational argument.
\qed
The proof of our main result on the parabolic-elliptic limit in (\ref{0eps}) is now almost immediate:\abs
\proofc of Theorem \ref{theo47}. \quad
  With $(\eps_{j_k})_{k\in\N}$ and $(n,c,u,P)$ taken as in Lemma \ref{lem38} and Lemma \ref{lem39},
  from Lemma \ref{lem38} and Lemma \ref{lem45} we immediately infer that (\ref{47.3}), (\ref{47.4}) and (\ref{47.7})  hold, and that moreover
  \be{47.99}
	\ceps \to c
	\qquad \mbox{in } L^\infty_{loc}((0,T];L^2(\Omega)) \cap L^2_{loc}((0,T];W^{1,2}(\Omega))
  \ee
  as $\eps=\eps_{j_k}\searrow 0$.
  Observing that for each fixed $\hatq>N$, Lemma \ref{lem50} provides $C_1>0$ such that
  \be{47.999}
	\|\ceps(\cdot,t)\|_{W^{1,\hatq}(\Omega)}
	\le C_1
	\qquad \mbox{for all $t\in (0,T)$ and each } \eps\in (\eps_j)_{j\in\N},
  \ee
  we firstly obtain (\ref{47.6}) as an immediate consequence thereof, and can furthermore secondly complete the derivation of
  (\ref{47.5}): As $\hatq>N$, namely, the Gagliardo-Nirenberg inequality yields $a\in (0,1)$ and $C_2>0$ such that
  for all $t\in (0,T)$ and $\eps\in (\eps_j)_{j\in\N}$,
  \bas
	\|\ceps(\cdot,t)-c(\cdot,t)\|_{C^0(\bom)}
	&\le& C_2 \|\ceps(\cdot,t)-c(\cdot,t)\|_{W^{1,\hatq}(\Omega)}^a \|\ceps(\cdot,t)-c(\cdot,t)\|_{L^2(\Omega)}^{1-a} \\
	&\le& C_2 \Big\{ \|\ceps(\cdot,t)\|_{W^{1,\hatq}(\Omega)} + \|c(\cdot,t)\|_{W^{1,\hatq}(\Omega)} \Big\}^a
	\cdot \|\ceps(\cdot,t)-c(\cdot,t)\|_{L^2(\Omega)}^{1-a},
  \eas
  so that combining (\ref{47.99}) with (\ref{47.999}) shows that indeed $\ceps\to c$ also
  in $L^\infty_{loc}((0,T);C^0(\bom))$.
  Finally, Lemma \ref{lem39} in conjunction with Lemma \ref{lem42} and Lemma \ref{lem46} guarantees that in fact
  $(n,c,u,P)$ solves (\ref{0e}) classically in $\Omega\times (0,T)$.
\qed
\mysection{Small-data solutions to an unforced Keller-Segel-Navier-Stokes system.
Proof of Theorem \ref{theo62}}\label{sect_t62}
The purpose of this section consists in providing a first exemplary application of Theorem \ref{theo47}, namely
in the framework of Theorem \ref{theo62}.
To this end, as for general $S$ the no-flux boundary conditions in (\ref{0eps}) need not reduce to separate homogeneous
Neumann boundary conditions for $\neps$ and $\ceps$, following \cite{win_ct_rot_SIMA} we introduce an appropriate regularization in which $S$ vanishes near the
lateral boundary. More precisely, let us fix $(\rho_\eta)_{\eta\in
(0,1)} \subset C_0^\infty(\Omega)$ and $(\chi_\eta)_{\eta\in
(0,1)} \subset C^\infty([0,\infty))$ such that \bas
    0 \le \rho_\eta \le 1 \mbox{ in $\Omega$ \quad with \quad}
    \rho_\eta \nearrow 1 \mbox{ in $\Omega$ as $\eta\searrow 0$},
\eas
and that
\be{chie}
    0 \le \chi_\eta \le 1 \mbox{ in $[0,\infty)$ \quad with \quad}
    \chi_\eta\equiv 0 \mbox{ in $[\frac{1}{\eta},\infty)$\, and } \,\,
    \chi_\eta \nearrow 1 \mbox{ in $[0,\infty)$ as $\eta\searrow 0$.}
\ee
For $\eta\in (0,1)$, we then define
\be{Seps}
    S_\eta(x,n,c):=\rho_\eta(x) \cdot \chi_\eta(n) \cdot S(x,n,c),
    \qquad (x,n,c)\in \bar\Omega\times [0,\infty)^2,
\ee and observe that $S_\eta\in C^2(\bar\Omega\times
[0,\infty)^2;\R^{N\times N})$.    \abs
Given $\eps>0$, for $\eta\in (0,1)$ we consider the approximate versions of (\ref{0eps}) given by
\be{0epseta}
	\left\{ \begin{array}{lcll}
     \, \, 	\partial_t \nee + \uee\cdot\nabla \nee
	&=& \Delta \nee - \nabla \cdot \Big( \nee S_\eta(x,\nee,\cee)\cdot\nabla \cee\Big),
	\qquad & x\in \Omega, \ t>0, \\[1mm]
	\eps \partial_t \cee + \uee\cdot\nabla \cee
	&=& \Delta \cee - \cee + \nee,
	\qquad & x\in \Omega, \ t>0, \\[1mm]
	\partial_t \uee + \kappa (\uee\cdot\nabla)\uee
	&=&  \Delta \uee + \nabla \Pee + \nee \nabla\phi,
	\qquad \nabla\cdot \uee=0,
	\qquad & x\in \Omega, \ t>0, \\[1mm]
	& & \hspace*{-46mm}
	\frac{\partial\nee}{\partial\nu}=\frac{\partial\cee}{\partial\nu}=0,
	\quad \uee=0,
	\qquad & x\in\pO, \ t>0, \\[1mm]
	& & \hspace*{-45mm}
	\nee(x,0)=n_0(x), \quad \cee(x,0)=c_0(x), \quad \uee(x,0)=u_0(x),
	\qquad & x\in\Omega,
	\end{array} \right.
\ee
Using that $S_\eta(x,n,c)=0$ whenever $n\ge \frac{1}{\eta}$, by following a series of standard arguments
(see e.g.\cite{win_CPDE2012}) one can readily verify that for each $\eps>0$ and $\eta\in (0,1)$
this problem possesses a globally defined classical solution $(\nee,\cee,\uee,\Pee)$ for which
$\nee$ and $\cee$ are nonnegative in $\Omega\times (0,\infty)$.\abs
In order to derive appropriate bounds for these solutions, independently of $\eps$ and $\eta$, we start by again using
Lemma \ref{lem01} to refine the differential inequality appearing therein as follows.
\begin{lem}\label{lem1}
  Let $N\ge 2, p>N, q>2$ and $r>N$.
  Then there exists $K_1(p,q,r)>0$ such that for all $\eps>0$ and any $\eta\in (0,1)$,
  \bea{1.1}
	\frac{\eps}{q} \frac{d}{dt} \io |\nabla\cee|^q
	&+& \Big\{ \frac{1}{4} - K_1(p,q,r) \|\uee\|_{L^r(\Omega)}^2 \Big\} \cdot \io |\nabla\cee|^{q-2} |D^2 \cee|^2 \nn\\
	&+& \Big\{ 1 - \frac{1}{q^2} - K_1(p,q,r) \|\uee\|_{L^r(\Omega)}^2 \Big\} \cdot \io |\nabla\cee|^q \nn\\[2mm]
	&\le& K_1(p,q,r) \|\nee\|_{L^p(\Omega)}^q
	\qquad \mbox{for all } t>0.
  \eea
\end{lem}
\proof
  Using Lemma \ref{lem01} as a starting point, thanks to  (\ref{01.1})  we can fix $C_1=C_1(p,q,r)>0$ such that for all
  $\eps>0$ and $\eta\in (0,1)$,
  \bea{1.2}
	& & \hspace*{-30mm}
	\frac{\eps}{q} \frac{d}{dt} \io |\nabla\cee|^q
	+ \frac{1}{4} \io |\nabla\cee|^{q-2} |D^2\cee|^2
	+ \Big(1-\frac{1}{q^2} \Big) \io |\nabla \cee|^q \nn\\
	&\le& C_1 \|\nee\|_{L^p(\Omega)}^q
	+ C_1 \|\uee\|_{L^r(\Omega)}^2 \|\nabla\cee\|_{L^\frac{qr}{r-2}(\Omega)}^q \\
	&=& C_1 \|\nee\|_{L^p(\Omega)}^q
	+ C_1 \|\uee\|_{L^r(\Omega)}^2 \Big\| |\nabla\cee|^\frac{q}{2} \Big\|_{L^\frac{2r}{r-2}(\Omega)}^2
	\qquad \mbox{for all } t>0.
  \eea
  Here since $r>N$ and thus $\frac{2r}{r-2}<\frac{2N}{N-2}$, we may use the continuity of the embedding
  $W^{1,2}(\Omega)\hra L^\frac{2r}{r-2}(\Omega)$ to find $C_2=C_2(p,q,r)>0$ such that
  \bas
	C_1 \|\uee\|_{L^r(\Omega)}^2 \Big\| |\nabla\cee|^\frac{q}{2} \Big\|_{L^\frac{2r}{r-2}(\Omega)}^2
	&\le& C_2 \|\uee\|_{L^r(\Omega)}^2
	\cdot \bigg\{ \Big\| \nabla |\nabla\cee|^\frac{q}{2} \Big\|_{L^2(\Omega)}^2
	+ \Big\| |\nabla\cee|^\frac{q}{2} \Big\|_{L^2(\Omega)}^2 \bigg\} \\
	&=& C_2\|\uee\|_{L^r(\Omega)}^2
	\cdot \bigg\{ \frac{q^2}{4} \io |\nabla\cee|^{q-4} |D^2 \cee \cdot \nabla \cee|^2
	+ \io |\nabla\cee|^q \bigg\} \\
	&\le& \frac{q^2}{4} C_2 \|\uee\|_{L^r(\Omega)}^2 \io |\nabla\cee|^{q-2} |D^2 \cee|^2
	+ C_2 \|\uee\|_{L^r(\Omega)}^2 \io |\nabla \cee|^q
  \eas
  for all $t>0$.
  When inserted into (\ref{1.2}), this yields (\ref{1.1}) on letting $K_1(p,q,r):=\max \big\{C_1, \frac{q^2}{4}C_2 \big\}$, for
  instance.
\qed
In consequence, if $\uee$ is suitably small, then also $\nabla\cee$ can be estimated in a favorable manner:
\begin{lem}\label{lem2}
  Let $N\ge 2, p>N$, $q>2$ and $r>N$.
  Then there exist $\delta_1(p,q,r)>0$ and $K_2(p,q,r)>0$ such that if $\eps>0$, $\eta\in (0,1)$ and $T>0$ are such that
  \be{2.1}
	\|\uee(\cdot,t)\|_{L^r(\Omega)} \le \delta_1(p,q,r)
	\qquad \mbox{for all } t\in (0,T),
  \ee
  then
  \be{2.2}
	\|\nabla\cee(\cdot,t)\|_{L^q(\Omega)}
	\le \max \bigg\{ \|\nabla c_0\|_{L^q(\Omega)} \, , \,
	K_2 \cdot \sup_{s\in (0,t)} \|\nee(\cdot,s)\|_{L^p(\Omega)} \bigg\}
	\qquad \mbox{for all } t\in (0,T).
  \ee
\end{lem}
\proof
  With $K_1(p,q,r)>0$ taken from Lemma \ref{lem1}, we let $\delta_1(p,q,r):=\frac{1}{\sqrt{4K_1(p,q,r)}}$ and then obtain
  from (\ref{1.1}) that if (\ref{2.1}) holds for some $\eps>0$, $\eta\in (0,1)$ and $T>0$, then
  $y(t):=\io |\nabla\cee(\cdot,t)|^q$, $t\ge 0$, satisfies
  \bas
	\frac{\eps}{q} \cdot y'(t) + \frac{1}{2} y(t) \le K_1(p,q,r) \|\neps(\cdot,t)\|_{L^p(\Omega)}^q
	\qquad \mbox{for all } t\in (0,T).
  \eas
  Therefore,
  if given any $t\in (0,T)$ we let $M(t):=K_1(p,q,r) \cdot \sup_{s\in (0,t)} \|\neps(\cdot,s)\|_{L^p(\Omega)}^q$, then
  \bas
	\frac{\eps}{q} \cdot y'(s) + \frac{1}{2} y(s) \le M(t)
	\qquad \mbox{for all } s\in (0,t),
  \eas
  so that a comparison argument yields the inequality
  \bas
	y(s) \le \max \Big\{ y(0) \, , \, 2M(t) \Big\}
	\qquad \mbox{for all } s\in [0,t].
  \eas
  When evaluated at $s=t$, this precisely leads to (\ref{2.2}) upon defining $K_2(p,q,r):=(2K_1(p,q,r))^\frac{1}{q}$.
\qed
Now the above hypothesis can be fulfilled if $u_0$ and $\nee$ are appropriately small:
\begin{lem}\label{lem3}
  Let $N\ge 2, p>1 $ and $r>N$ be such that
  \be{3.1}
	p>\frac{Nr}{N+2r}.
  \ee
  Then for all $\delta>0$ there exists $\delta_3(\delta,p,r)>0$ such that if
  \be{3.2}
	\|u_0\|_{L^r(\Omega)} \le \delta_3(\delta,p,r),
  \ee
  and if for some $\eps>0, \eta\in (0,1)$ and $T>0$ we have
  \be{3.3}
	\|\nee(\cdot,t)\|_{L^p(\Omega)} \le \delta_3(\delta,p,r)
	\qquad \mbox{for all } t\in (0,T),
  \ee
  then
  \be{3.4}
	\|\uee(\cdot,t)\|_{L^r(\Omega)} \le \delta	
	\qquad \mbox{for all } t\in (0,T).
  \ee
\end{lem}
\proof
  By relying on known regularization features of the Stokes semigroup $(e^{-tA})_{t\ge 0}$ (\cite{giga1986}),
  let us fix $C_1=C_1(r)>0$, $C_2=C_2(r)>0$, $C_3=C_3(p,r)>0$ and $\mu>0$ such that for all $t>0$,
  \be{3.44}
	\|e^{-tA} \varphi\|_{L^r(\Omega)} \le C_1\|\varphi\|_{L^r(\Omega)}
	\qquad \mbox{for all }   \varphi \in L^r_\sigma(\Omega)
  \ee
  and
  \be{3.5}
	\|e^{-tA} \proj \nabla \cdot\varphi \|_{L^r(\Omega)}
	\le C_2 t^{-\frac{1}{2} - \frac{N}{2r}} e^{-\mu t} \|\varphi\|_{L^\frac{r}{2}(\Omega)}
	\qquad \mbox{for all $\varphi\in C^1(\bom;\R^{N\times N})$ such that $\varphi=0$ on $\pO$}
  \ee
  as well as
  \be{3.6}
	\|e^{-tA} \proj \varphi\|_{L^r(\Omega)}
	\le C_3 t^{-\frac{N}{2}(\frac{1}{p}-\frac{1}{r})_+} e^{-\mu t} \|\varphi\|_{L^p(\Omega)}
	\qquad \mbox{for all } \varphi\in C^0(\bom; \R^N).
  \ee
  We then fix $\delta>0$ and may without loss of generality assume that
  \be{3.7}
	C_2 C_4 |\kappa| \delta^2 \le \frac{\delta}{6},
  \ee
  where $C_4=C_4(r):=\int_0^\infty \sigma^{-\frac{1}{2}-\frac{N}{2r}} e^{-\mu \sigma} d\sigma <\infty$ since $r>N$.
  Noting that thanks to (\ref{3.1}) we moreover know that also
  $C_5=C_5(p,r):=\int_0^\infty \sigma^{-\frac{N}{2}(\frac{1}{p}-\frac{1}{r})_+} e^{-\mu\sigma} d\sigma$ is finite,
  we thereupon pick $\delta_3=\delta_3(\delta,p,r)>0$ small enough such that both
  \be{3.8}
	C_1 \delta_3 \le \frac{\delta}{6}
  \ee
  and
  \be{3.9}
	C_3 C_5 \|\nabla\phi\|_{L^\infty(\Omega)} \delta_3 \le \frac{\delta}{6}
  \ee
  hold, and suppose that (\ref{3.2}) and (\ref{3.3}) are satisfied with some $\eps>0, \eta\in (0,1)$ and $T>0$.
  Then since $\uee$ clearly is a mild solution of its respective subproblem in (\ref{0epseta}), we may use (\ref{3.44}),
  (\ref{3.5}) and (\ref{3.6}) along with the H\"older inequality to estimate
  \bea{3.10}
	& & \hspace*{-16mm}
	\|\uee(\cdot,t)\|_{L^r(\Omega)} \nn\\
	&=& \Bigg\| e^{-tA} u_0
	- \kappa \int_0^t e^{-(t-s)A} \proj \nabla \cdot ( \uee(\cdot,s) \mult \uee(\cdot,s)) ds
	+ \int_0^t e^{-(t-s)A} \proj [n(\cdot,s)\nabla \phi] ds \Bigg\|_{L^r(\Omega)} \nn\\
	&\le& C_1 \|u_0\|_{L^r(\Omega)}
	+ C_2 |\kappa| \int_0^t (t-s)^{-\frac{1}{2}-\frac{N}{2r}} e^{-\mu(t-s)}
	\|\uee(\cdot,s) \mult \uee(\cdot,s)\|_{L^\frac{r}{2}(\Omega)} ds \nn\\
	& & + C_3 \int_0^t (t-s)^{-\frac{N}{2}(\frac{1}{p}-\frac{1}{r})_+} e^{-\mu(t-s)}
		\|\nee(\cdot,s) \nabla\phi\|_{L^p(\Omega)} ds \nn\\
	&\le& C_1\|u_0\|_{L^r(\Omega)}
	+ C_2|\kappa| \int_0^t (t-s)^{-\frac{1}{2} - \frac{N}{2r}} e^{-\mu(t-s)} \|\uee(\cdot,s)\|_{L^r(\Omega)}^2 ds \nn\\
	& & + C_3\|\nabla \phi\|_{L^\infty(\Omega)} \int_0^t (t-s)^{-\frac{N}{2}(\frac{1}{r}-\frac{1}{p})_+} e^{-\mu(t-s)}
		\|\nee(\cdot,s)\|_{L^p(\Omega)} ds
	\qquad \mbox{for all } t\in [0,T].
  \eea
  In order to verify that this implies the inequality
  \be{3.11}
	M(T_0) <\delta
	\qquad \mbox{for all } T_0\in [0,T]
  \ee
  for
  \bas
	M(T_0):=\sup_{t\in [0,T_0]} \|\uee(\cdot,t)\|_{L^r(\Omega)},
	\qquad T_0\in [0,T],
  \eas
  assuming (\ref{3.11}) to be false we could make use of the continuity of $M$ and e.g.~combine (\ref{3.10}) with (\ref{3.7})
  to find $T_\star\in (0,T]$ such that $M(T_0)<\delta$ for all $T_0\in [0,T_\star)$ but $M(T_\star)=\delta$.
  According to (\ref{3.10}) and our hypotheses (\ref{3.2}) and (\ref{3.3}) in conjunction with (\ref{3.7}), (\ref{3.8})
  and (\ref{3.9}), however, this would mean that
  \bas
	\delta=M(T_\star)
	&\le& C_1 \delta_3 + C_2|\kappa| \delta^2 \int_0^t (t-s)^{-\frac{1}{2}-\frac{N}{2r}} e^{-\mu(t-s)} ds \\
	& & + C_3 \|\nabla\phi\|_{L^\infty(\Omega)} \delta_3 \int_0^t
		(t-s)^{-\frac{N}{2}(\frac{1}{p}-\frac{1}{r})_+} e^{-\mu(t-s)} ds \\
	&\le& C_1 \delta_3 + C_2 C_4 |\kappa| \delta^2
	+ C_3 C_5 \|\nabla \phi\|_{L^\infty(\Omega)} \delta_3 \\
	&\le& \frac{\delta}{6} + \frac{\delta}{6} + \frac{\delta}{6}=\frac{\delta}{2},
  \eas
  which is absurd. Therefore, (\ref{3.11}) and hence (\ref{3.4}) must be valid.
\qed
The hypotheses of the latter lemma, however, are satisfied if $\nabla\cee$ and $\uee$ are conveniently small:
\begin{lem}\label{lem4}
  Let $N\ge 2$. Then for all $p>1, q>N$ and $r>N$ there exists $\delta_2(p,q,r)>0$ such that if for some $\eps>0$,
  $\eta\in (0,1)$ and $T>0$ we have
  \be{4.1}
	\|\nabla \cee(\cdot,t)\|_{L^q(\Omega)} \le \delta_2(p,q,r)
	\qquad \mbox{for all } t\in (0,T)
  \ee
  and
  \be{4.2}
	\|\uee(\cdot,t)\|_{L^r(\Omega)} \le \delta_2(p,q,r)
	\qquad \mbox{for all } t\in (0,T),
  \ee
  then
  \be{4.4}
	\|\nee(\cdot,t)\|_{L^p(\Omega)} \le 2\|n_0\|_{L^p(\Omega)}
	\qquad \mbox{for all } t\in (0,T).
  \ee
\end{lem}
\proof
  Given $p>1, q>N$ and $r>N$, by using a well-known smoothing property of the Neumann heat semigroup
  $(e^{t\Delta})_{t\ge 0}$
  (\cite{win_JDE2010}) we can fix $\mu>0, C_1=C_1(p,q)>0$ and $C_2=C_2(p,r)>0$ such that whenever $t>0$,
  \bea{4.5}
	\|e^{t\Delta} \nabla\cdot\varphi\|_{L^p(\Omega)}
	&\le& \min \Big\{ C_1 t^{-\frac{1}{2}-\frac{N}{2q}} e^{-\mu t} \|\varphi\|_{L^\frac{pq}{p+q}(\Omega)} \, , \,
	C_2 t^{-\frac{1}{2}-\frac{N}{2r}} e^{-\mu t} \|\varphi\|_{L^\frac{pr}{p+r}(\Omega)} \Big\} \nn\\
	& & \hspace*{10mm}
	\qquad \mbox{for all $\varphi\in C^1(\bom;\R^N)$ such that $\varphi\cdot\nu=0$ on } \pO,
   \eea
  and thereupon let $\delta_2=\delta_2(p,q,r)>0$ be small enough such that
  \be{4.6}
	C_1 C_3 K_S \delta_2 \le \frac{1}{4}
  \ee
  and
  \be{4.7}
	C_2 C_4 \delta_2 \le \frac{1}{4},
  \ee
  where $C_2=C_2(q):=\int_0^\infty \sigma^{-\frac{1}{2}-\frac{N}{2q}} e^{-\mu \sigma} d\sigma$ and
  where $C_3=C_3(r):=\int_0^\infty \sigma^{-\frac{1}{2}-\frac{N}{2r}} e^{-\mu \sigma} d\sigma$ are finite due to our
  assumptions that $q>N$ and $r>N$.\abs
  Then assuming (\ref{4.1}) and (\ref{4.2}) to be valid for some $\eps>0, \eta\in (0,1)$ ad $T>0$, we may employ a
  variation-of-constants representation associated with the first equation in (\ref{0epseta}) to see that thanks to the contractivity of $(e^{t\Delta})_{t\ge 0}$ on $L^p(\Omega)$,
  \bea{4.8}
	& & \|\nee(\cdot,t)\|_{L^p(\Omega)}   \nn \\
	&=& \Bigg\| e^{t\Delta} n_0
	- \int_0^t e^{(t-s)\Delta}
	\nabla\cdot \Big( \nee(\cdot,s) S_\eta(\cdot,\nee(\cdot,s),\cee(\cdot,s)\big) \cdot \nabla \cee(\cdot,s) \Big) ds \nn\\
	& & \hspace*{10mm}
	- \int_0^t e^{(t-s)\Delta} \nabla \cdot \Big( \nee(\cdot,s) \uee(\cdot,s)\Big) ds \Bigg\|_{L^p(\Omega)} \nn\\
	&\le& \|n_0\|_{L^p(\Omega)}
	+ C_1 \int_0^t (t-s)^{-\frac{1}{2}-\frac{N}{2q}} e^{-\mu(t-s)}
	\Big\| \nee(\cdot,s) S_\eta(\cdot,\nee(\cdot,s),\cee(\cdot,s)) \cdot \nabla \cee(\cdot,s) \Big\|_{L^\frac{pq}{p+q}(\Omega)} ds \nn\\
	& & + C_2 \int_0^t (t-s)^{-\frac{1}{2}-\frac{N}{2r}} e^{-\mu(t-s)}
	\|\nee(\cdot,s)\uee(\cdot,s)\|_{L^\frac{pr}{p+r}(\Omega)} ds
  \eea
  for all $t\in (0,T)$.
  Here by the H\"older inequality,   (\ref{S}), from (\ref{4.1}) we know that abbreviating
  $M:=\|\nee\|_{L^\infty((0,T);L^p(\Omega))}$ we have
  \bas
	\Big\|\nee(\cdot,s) S_\eta(\cdot,\nee(\cdot,s),\cee(\cdot,s)) \cdot \nabla \cee(\cdot,s)
		\Big\|_{L^\frac{pq}{p+q}(\Omega)}
	&\le& K_S \|\nee(\cdot,s)\|_{L^p(\Omega)} \|\nabla\cee(\cdot,s)\|_{L^q(\Omega)} \\
	&\le& K_S \delta_2 M
	\qquad \mbox{for all } s\in (0,T),
  \eas
  while similarly (\ref{4.2}) guarantees that
  \bas
	\|\nee(\cdot,s)\uee(\cdot,s)\|_{L^\frac{pr}{p+r}(\Omega)}
	&\le& \|\nee(\cdot,s)\|_{L^p(\Omega)} \|\uee(\cdot,s)\|_{L^r(\Omega)} \\
	&\le& \delta_2 M
	\qquad \mbox{for all } s\in (0,T).
  \eas
  Therefore, we may use (\ref{4.6}) and (\ref{4.7}) to infer from (\ref{4.8}) that
  \bas
	\|\nee(\cdot,t)\|_{L^p(\Omega)}
	&\le& \|n_0\|_{L^p(\Omega)}
	+ C_1 K_S \delta_2 M \int_0^t (t-s)^{-\frac{1}{2}-\frac{N}{2q}} e^{-\mu(t-s)} ds \\
	& & + C_2 \delta_2 M \int_0^t (t-s)^{-\frac{1}{2}-\frac{N}{2r}} e^{-\mu(t-s)} ds \\
	&\le& \|n_0\|_{L^p(\Omega)}	
	+ C_1 C_3 K_S \delta_2 M
	+ C_2 C_4 \delta_2 M \\
	&\le& \|n_0\|_{L^p(\Omega)}
	+ \frac{M}{4} + \frac{M}{4}
	\qquad \mbox{for all } t\in (0,T),
  \eas
  which implies that
  \bas
	M \le \|n_0\|_{L^p(\Omega)} + \frac{M}{2}
  \eas
  and hence completes the proof.
\qed
Now a self-map type argument combines the latter lemmata so as to make sure that for suitably small initial data,
all the above assumptions can be fulfilled simultaneously:
\begin{lem}\label{lem61}
  Let $N\ge 2, p>\max\{2,N\}, q>N$ and $r>N$.
  Then there exists $C=C(p,q,r)>0$ such that if $n_0, c_0$ and $u_0$ satisfy (\ref{init}) with
  \be{61.1}
	\|n_0\|_{L^p(\Omega)} \le \frac{1}{C},
	\qquad
	\|\nabla c_0\|_{L^q(\Omega)} \le \frac{1}{C}
	\qquad \mbox{and} \qquad
	\|u_0\|_{L^r(\Omega)} \le \frac{1}{C},
  \ee
  then for all $\eps>0$ and $\eta\in (0,1)$, the solution of (\ref{0epseta}) has the properties that
  \be{61.2}
	\|\nee(\cdot,t)\|_{L^p(\Omega)} \le C,
	\quad  \, \,
	\|\nabla \cee(\cdot,t)\|_{L^q(\Omega)} \le C
	\quad  \, \,  \mbox{and} \quad  \,  \,
	\|\uee(\cdot,t)\|_{L^r(\Omega)} \le C
	\quad  \,  \,   \mbox{for all } \,  t>0.
  \ee
\end{lem}
\proof
  Given $p>\max\{2,N\}, q>N$ and $r>N$, we take  $\delta_1=\delta_1(p,q, r)>0$ and $K_2=K_2(p,q,r)>0$ from Lemma \ref{lem2}
  and let $\delta_2=\delta_2(p,q,r)>0$ be as provided by Lemma \ref{lem4}. Then since $p>N>\frac{Nr}{N+2r}$,
  an application of Lemma \ref{lem3} to $\delta:=\min\{\delta_1,\delta_2\}$ yields $\delta_3=\delta_3(p,q,r)>0$
  with the property that whenever (\ref{3.2}) and (\ref{3.3}) hold for some $\eps>0, \eta\in (0,1)$ and $T>0$, we have
  \be{61.3}
	\|\uee(\cdot,t)\|_{L^r(\Omega)} \le \delta_1
	\qquad \mbox{for all } t\in (0,T)
  \ee
  and
  \be{61.4}
	\|\uee(\cdot,t)\|_{L^r(\Omega)} \le \delta_2
	\qquad \mbox{for all } t\in (0,T).
  \ee
  We now suppose that $n_0,c_0$ and $u_0$ comply with (\ref{init}) and are such that
  \be{61.5}
	3K_2 \|n_0\|_{L^p(\Omega)} \le \delta_2
  \ee
  and
  \be{61.6}
	\|\nabla c_0\|_{L^q(\Omega)} \le \delta_2
  \ee
  as well as
  \be{61.7}
	\|u_0\|_{L^r(\Omega)} \le \delta_3
  \ee
  and
  \be{61.77}
	3\|n_0\|_{L^p(\Omega)} \le \delta_3,
  \ee
  and we claim that then for each $\eps>0$ and $\eta\in (0,1)$, the obviously well-defined element
  \bas
	T\equiv T_ {\eps \eta}
	:=\sup \bigg\{ \widehat{T}>0 \ \bigg| \ \|\nee(\cdot,t)\|_{L^p(\Omega)} < 3\|n_0\|_{L^p(\Omega)}
	\ \mbox{for all } t\in (0,\widehat{T}) \bigg\}
  \eas
  of $(0,\infty]$ actually satisfies $T_{\eps \eta}=\infty$.\abs
  To see this, we note that by definition of $T$,
  \be{61.8}
	\|\nee(\cdot,t)\|_{L^p(\Omega)} <  3\|n_0\|_{L^p(\Omega)}
	\qquad \mbox{for all } t\in (0,T),
  \ee
  which in conjunction with (\ref{61.7}) and (\ref{61.77}) allows for an application of Lemma \ref{lem3} to conclude
  that in fact both (\ref{61.3}) and (\ref{61.4}) hold. In particular, (\ref{61.3}) enables us to employ Lemma \ref{lem2}
  to see that thanks to (\ref{61.6}), again (\ref{61.8}), and (\ref{61.5}),
  \bea{61.9}
	\|\nabla \cee(\cdot,t)\|_{L^q(\Omega)}
	&\le& \max \Big\{ \|\nabla c_0\|_{L^q(\Omega)} \, , \, K_2 \cdot 3\|n_0\|_{L^p(\Omega)} \Big\} \nn\\[2mm]
	&\le& \delta_2
	\qquad \mbox{for all } t\in (0,T),
  \eea
  which in turn, when combined with (\ref{61.4}), makes it possible to infer from Lemma \ref{lem4} that
  \bas
	\|\nee(\cdot,t)\|_{L^p(\Omega)} \le 2\|n_0\|_{L^p(\Omega)}
	\qquad \mbox{for all } t\in (0,T).
  \eas
  As $n_0\not\equiv 0$ by (\ref{init}), by continuity of $\nee$ this shows that indeed $T_{\eps \eta}$ cannot be
  finite for any $\eps>0$ and $\eta\in (0,1)$, and that thus (\ref{61.2}) results as a consequence of (\ref{61.8}),
  (\ref{61.9}) and (\ref{61.3}) if in accordance with (\ref{61.5})-(\ref{61.77}), the constant $C$ in (\ref{61.1})
  and (\ref{61.2}) is chosen suitably large.
\qed
In fact, we have thereby proved the essential body of Theorem \ref{theo62} already:\abs
\proofc of Theorem \ref{theo62}. \quad
  According to Lemma \ref{lem61}, there exists $\delta=\delta(p,q,r)>0$ such that (\ref{62.1}) implies the boundedness
  properties in (\ref{61.2}) uniformly with respect to $\eps>0$ and $\eta\in (0,1)$.
  Thanks to the estimates  thereby implied through  Lemma \ref{lem32}, Lemma \ref{lem33} and Lemma \ref{lem45},
  by means of a standard subsequence extraction procedure this can readily be seen to entail,
  for each $\eps>0$, the existence of a global classical solution $(\neps,\ceps,\ueps,\Peps)$ of (\ref{0eps}) which
  in fact has the properties that $n_{\eps \eta_l} \to \neps$, $c_{\eps \eta_l} \to \ceps$ and $u_{\eps \eta_l} \to \ueps$
  a.e.~in $\Omega \times (0,\infty)$ with some $(\eta_l)_{l\in\N} \subset (0,1)$ such that $\eta_l\searrow 0$
  as $l\to\infty$ (\cite{cao_lankeit}).\\
  The remaining part of the statement then directly results from Theorem \ref{theo47} and the boundedness features of
  $(\nabla \ceps)_{\eps>0}$ and $(\ueps)_{\eps>0}$ implied by (\ref{61.2}).
\qed
\mysection{A logistic Keller-Segel system. Proof of Theorems \ref{theo66} and \ref{theo67}}\label{sect8}
As a second example for taking a parabolic-elliptic limit along the lines of Theorem \ref{theo47}, in this section
we shall consider the one-dimensional logistic Keller-Segel system (\ref{L}) for fixed $D>0, a\in\R, b\ge 0$ and $\eps>0$.\abs
Again we start by stating an almost immediate basic property.
\begin{lem}\label{lem63}
  Let $T>0$. Then there exists $C(T)>0$ such that for any $\eps>0$,
  \be{63.2}
	\int_0^1 \ceps(\cdot,t) \le C(T)
	\qquad \mbox{for all } t\in (0,T)
  \ee
\end{lem}
\proof
  As an immediate consequence of Lemma \ref{lem31}, we obtain $C_1(T)>0$ such that
  \bas
	\int_0^1 \neps (\cdot,t)  \le C_1(T)		
	\qquad \mbox{for all } t\in (0,T).
  \eas
  Thereupon, using (\ref{L}) we can estimate
  \bas
	\eps \frac{d}{dt} \int_0^1 \ceps (\cdot,t) + \int_0^1 \ceps (\cdot,t) = \int_0^1 \neps (\cdot,t) \le C_1(T)
	\qquad \mbox{for all } t\in (0,T),
  \eas
  which by comparison implies that
  \bas
	\int_0^1 \ceps (\cdot,t) \le \max \bigg\{ \int_0^1 c_0 \, , \, C_1(T) \bigg\}
	\qquad \mbox{for all } t\in (0,T),
  \eas
  as intended.
\qed
Now in the spatially one-dimensional setting considered here, the availability of favorable embeddings allows us to conclude the following from an essentially well-established testing procedure.
\begin{lem}\label{lem64}
  Let $T>0$. Then there exists $C(T)>0$ such that
  \be{64.1}
	\int_0^T \int_0^1 c_{\eps xx}^2 \le C(T)
	\qquad \mbox{for all } \eps \in (0,1).
  \ee
\end{lem}
\proof
By referring to both PDEs in (\ref{L}) and employing Young's inequality, we see that whenever $\eps>0$,
  \bea{64.2}
	& & \hspace*{-20mm}
	\frac{d}{dt} \bigg\{ \int_0^1 \neps \ln \neps (\cdot,t)  + \frac{\eps}{2}  \int_0^1 c_{\eps x}^2 (\cdot,t)  \bigg\}	
	+ D \int_0^1 \frac{n_{\eps x}^2}{\neps} (\cdot,t) + \int_0^1 c_{\eps xx}^2 (\cdot,t) + \int_0^1 c_{\eps x}^2 (\cdot,t)  \nn\\
	&=& - 2 \int_0^1 \neps c_{\eps xx} (\cdot,t)
	+ a \int_0^1 \neps \ln \neps (\cdot,t)
	- b \int_0^1 \neps^2 \ln \neps (\cdot,t)
	+ a \int_0^1  \neps (\cdot,t)
	- b \int_0^1  \neps^2 (\cdot,t)  \nn\\
	&\le& \frac{1}{2} \int_0^1 c_{\eps xx}^2 (\cdot,t)
	+ a \int_0^1 \neps \ln \neps (\cdot,t)
	- b \int_0^1 \neps^2 \ln \neps (\cdot,t)
	+ a \int_0^1 \neps (\cdot,t)
	+ \int_0^1 \neps^2 (\cdot,t)
  \eea
for all $t\in(0, T)$.   Since it can readily be verified by elementary analysis that thanks to the nonnegativity of $b$ there exists $C_1>0$
  with the property that
  \bas
	a\xi\ln \xi - b\xi^2 \ln \xi + a\xi + \xi^2
	\le 2\xi^2 + C_1
	\qquad \mbox{for all } \xi> 0,
  \eas
  and since the Gagliardo-Nirenberg inequality, Young's inequality and Lemma \ref{lem31} provide $C_2>0$ and $C_3(T)>0$
  such that for all $\eps>0$ we have
  \bas
	2 \int_0^1 \neps^2(\cdot,t)
	&=& 2\|\sqrt{\neps}\|_{L^4((0,1))}^4 \\
	&\le& C_2 \|(\sqrt{\neps})_x\|_{L^2((0,1))} \|\sqrt{\neps}\|_{L^2((0,1))}^3
	+ C_2 \|\sqrt{\neps}\|_{L^2((0,1))}^4 \\
	&\le& D \int_0^1 \frac{n_{\eps x}^2}{\neps}(\cdot,t)  + C_3(T)
	\qquad \mbox{for all } t\in (0,T),
  \eas
  from (\ref{64.2}) it thus follows that for any such $\eps$,
  \bas
	\frac{d}{dt} \bigg\{ \int_0^1 \neps \ln \neps (\cdot,t)  + \frac{\eps}{2} \int_0^1 c_{\eps x}^2 (\cdot,t)   \bigg\}	
	+ \frac{1}{2} \int_0^1 c_{\eps xx}^2 (\cdot,t)  \le C_4(T):=C_1 + C_3(T)
	\quad \mbox{for all } t\in (0,T).
  \eas
  Hence, when resorting to $\eps\in (0,1)$ we infer that
  \bas
	\int_0^1 \neps(\cdot,T)\ln \neps(\cdot,T)
	+ \frac{\eps}{2} \int_0^1 c_{\eps x}^2(\cdot,T)
	+ \frac{1}{2} \int_0^T \int_0^1 c_{\eps xx}^2
	&\le& \int_0^1 n_0\ln n_0 + \frac{\eps}{2} \int_0^1 c_{0x}^2 + C_4(T) T \\
	&\le& \int_0^1 n_0\ln n_0 + \frac{1}{2} \int_0^1 c_{0x}^2 + C_4(T) T,
  \eas
  which entails (\ref{64.1}) due to the fact that $\int_0^1 \neps(\cdot,T) \ln \neps(\cdot,T) \ge - \frac{1}{e}$.
\qed
In conjunction with the $L^1$ information from Lemma \ref{lem63}, the latter entails an estimate for $c_{\eps x}$
compatible with (\ref{47.22}):
\begin{lem}\label{lem65}
  For any $T>0$ one can find $C(T)>0$ with the property that
  \be{65.1}
	\int_0^T \|c_{\eps x}(\cdot,t)\|_{L^\infty((0,1))}^\frac{5}{2} dt \le C(T)
	\qquad \mbox{for all } \eps\in (0,1).
  \ee
\end{lem}
\proof
  As the Gagliardo-Nirenberg inequality says that with some $C_1>0$ we have
  \bas
	\|c_{\eps x}\|_{L^\infty((0,1))}^\frac{5}{2}
	\le C_1\|c_{\eps xx}\|_{L^2((0,1))}^2 \|\ceps\|_{L^1((0,1))}^\frac{1}{2}
	+ C_1\|\ceps\|_{L^1((0,1))}^\frac{5}{2}
	\qquad \mbox{for all $t>0$ and each } \eps>0,
  \eas
  the claim results upon integrating and combining Lemma \ref{lem64} with Lemma \ref{lem63}.
\qed
We can thereby directly pass to the limit $\eps\searrow 0$ by means of Theorem \ref{theo47}:\abs
\proofc of Theorem \ref{theo66}.\quad
  We pick any $q>5$ and then obtain as a particular consequence of Lemma \ref{lem65} that for each $T>0$,
  $(c_{\eps x})_{\eps\in (0,1)}$ is bounded in $L^\frac{5}{2}((0,T);L^q((0,1)))$.
  Since this choice of $q$ precisely ensures that $\frac{2}{5}+\frac{1}{2q}<\frac{1}{2}$, the conclusion follows by
  applying Theorem \ref{theo47} with $u_0\equiv 0$ and $\phi\equiv 0$, and recalling from standard literature
  (\cite{tello_win}, \cite{cieslak_win}) a well-known uniqueness property
  of (\ref{Le}) within the indicated class.
\qed
Thanks to a known result on spontaneous emergence of large densities in the limit problem (\ref{Le}) for suitably small
$D>0$, our statement from Theorem \ref{theo66} enables us to finally draw a similar conclusion
also for the fully parabolic problem when the parameter $\eps$ therein is appropriately small.\abs
\proofc of Theorem \ref{theo67}. \quad
  According to a result from \cite[Theorem 1.1]{win_JNLS} on the parabolic-elliptic problem
  (\ref{Le}), we can pick some nonnegative $n_0\in W^{1,\infty}((0,1))$ which is such that there exists $T>0$ having
  the property that to arbitrary $M>0$ there corresponds some $D_0>0$ such that for each $D\in (0,D_0)$, the solution
  $(n,c)\equiv (n_D,c_D)$ of (\ref{Le}) satisfies
  \be{67.2}
	n_D(x_0(D),t_0(D)) \ge 2M
  \ee
  with some $x_0(D)\in (0,1)$ and $t_0(D)\in (0,T)$.
  Now keeping $n_0, T$ and $M$ fixed, given any such $D$ and arbitrary nonnegative $c_0\in W^{1,\infty}((0,1))$
  we may employ Theorem \ref{theo66} to see that the associated solutions $(n_{D\eps}, c_{D\eps})$ of (\ref{L})
  approximate $(n_D,c_D)$ in the sense that, inter alia,   $n_{D\eps} \rightarrow n_D$ in $C^0([0,1]\times [0,T])$ as $\eps\searrow 0$.
  In particular, we can therefore find $\eps_0>0$ such that $n_{D\eps} \ge n_D-M$ in $(0,1)\times (0,T)$
  for all $\eps\in (0,\eps_0)$, which when evaluated at $(x_0(D),t_0(D))$ and combined with (\ref{67.2}) directly
  yields (\ref{67.1}).
\qed

\bigskip

{\bf Acknowledgements.} \quad
  Y. Wang was supported by the NNSF of China (no. 11501457) and Xihua University Scholars Training Program.  M.~Winkler acknowledges support of the {\em Deutsche Forschungsgemeinschaft} in the context of the project  {\em Analysis of chemotactic cross-diffusion in complex frameworks}.  Z. Xiang was partially supported by the NNSF of China under Grants 11571063 and 11771045.

\end{document}